\def\draft{n}
\documentclass[headings]{amsart}
\usepackage{fullpage,amssymb,epic,eepic,epsfig,amsbsy,amsmath,bbold}

\theoremstyle{plain}

\newtheorem{theorem}{Theorem}
\newtheorem{proposition}{Proposition}[section]
\newtheorem{lemma}[proposition]{Lemma}

\newtheorem{conjecture}{Conjecture}

\theoremstyle{definition}
\newtheorem{definition}[proposition]{Definition}

\theoremstyle{remark}

\newtheorem{remark}[proposition]{Remark}

\def\printname#1{
        \if\draft y
                \smash{\makebox[0pt]{\hspace{-0.5in}
                        \raisebox{8pt}{\tt\tiny #1}}}
        \fi
}

\newcommand{\psdraw}[2]
         {\begin{array}{c} \hspace{-1.3mm}
        \raisebox{-4pt}{\epsfig{figure=draws/#1.eps,width=#2}}
        \hspace{-1.9mm}\end{array}}

\newlength{\standardunitlength}
\catcode`\@=11
\long\def\@makecaption#1#2{%
     \vskip 10pt

\setbox\@tempboxa\hbox{
       \small\sf{\bfcaptionfont #1. }\ignorespaces #2}%
     \ifdim \wd\@tempboxa >\captionwidth {%
         \rightskip=\@captionmargin\leftskip=\@captionmargin
         \unhbox\@tempboxa\par}%
       \else
         \hbox to\hsize{\hfil\box\@tempboxa\hfil}%
     \fi}
\font\bfcaptionfont=cmssbx10 scaled \magstephalf
\newdimen\@captionmargin\@captionmargin=2\parindent
\newdimen\captionwidth\captionwidth=\hsize
\catcode`\@=12

\def\lbl#1{\label{#1}\printname{#1}}

\def\BN{\mathbb N}
\def\BZ{\mathbb Z}

\def\BQ{\mathbb Q}
\def\BR{\mathbb R}
\def\BC{\mathbb C}

\def\D{\Delta}

\def\P{\mathcal P}

\def\a{\alpha}

\def\l{\lambda}

\def\s{\sigma}

\def\ga{\gamma}

\def\w{\omega}

\def\e{\epsilon}

\def\th{\theta}
\def\Th{\Theta}
\def\s{\sigma}

\def\longto{\longrightarrow}
\def\w{\omega}
\def\tpsi{\tilde{\psi}}
\def\hpsi{\hat{\psi}}
\def\SL{\mathrm{SL}}
\def\pt{\partial}

\def\ch{\mathrm{ch}}
\def\xpe{x_p^{\epsilon}}

\begin{document}


\title[Asymptotics of $q$-difference equations]{
Asymptotics of $q$-difference equations}

\author{Stavros Garoufalidis}
\address{School of Mathematics \\
         Georgia Institute of Technology \\
         Atlanta, GA 30332-0160, USA}
\email{stavros@math.gatech.edu, 
URL: {\tt http://www.math.gatech.edu/$\sim$stavros }}

\author{Jeffrey S. Geronimo}
\address{School of Mathematics \\
         Georgia Institute of Technology \\
         Atlanta, GA 30332-0160, USA}
\email{geronimo@math.gatech.edu}

\thanks{The authors were supported in part by the National Science 
Foundation. \\
\newline
1991 {\em Mathematics Classification.} Primary 57N10. Secondary 57M25.
}

\date{April 20, 2004} 


\begin{abstract}
In this paper we develop an asymptotic analysis for formal and actual
solutions of $q$-difference equations, under a regularity assumption,
namely the non-collision and non-vanishing of the eigenvalues.
In particular, evaluations of regular solutions of regular 
$q$-difference equations 
have an exponential growth rate which can be computed from the $q$-difference
equation.

The motivation for the paper comes from a problem in Quantum Topology,
the Hyperbolic Volume Conjecture, which states that a sequence on Laurent
polynomials (the so-called colored Jones function of a knot), appropriately
evaluated, becomes a sequence of complex numbers that grows exponentially.
Moreover, the exponential growth rate is proportional to the volume of the
knot complement.

The connection of the Hyperbolic Volume Conjecture with the paper comes 
from the fact that the 
colored Jones function of a knot is a solution of a $q$-difference
equation, as was proven by TTQ. Le and the author.
\end{abstract}

\maketitle

\tableofcontents


\section{Introduction}
\lbl{sec.intro}

\subsection{The goal}
\lbl{sub.goal}

The goal of the paper is to intiate an approach to the Hyperbolic Volume 
Conjecture, via asymptotics of solutions of difference equations with a small
parameter. The Generalized Volume Conjecture links (conjecturally)
the (colored) Jones polynomial of a knot to hyperbolic geometry of its 
complement.

Since the colored Jones polynomial is a specific solution to a linear
$q$-difference equation, it follows that
the generalized volume conjecture is the WKB limit of a specific solution
of a linear difference equation with a small parameter.

Motivated by this, we study WKB asymptotics of formal and actual solutions
of difference equations with a small parameter, under certain regularity 
asymptions.

\subsection{The colored Jones function}
\lbl{sub.jones}

A knot in 3-space is a smooth embedding of a circle, considered up to
isotopy. Two of the simplest knots, the Trefoil ($3_1$) and the Figure Eight
($4_1$) are shown here:
$$
\psdraw{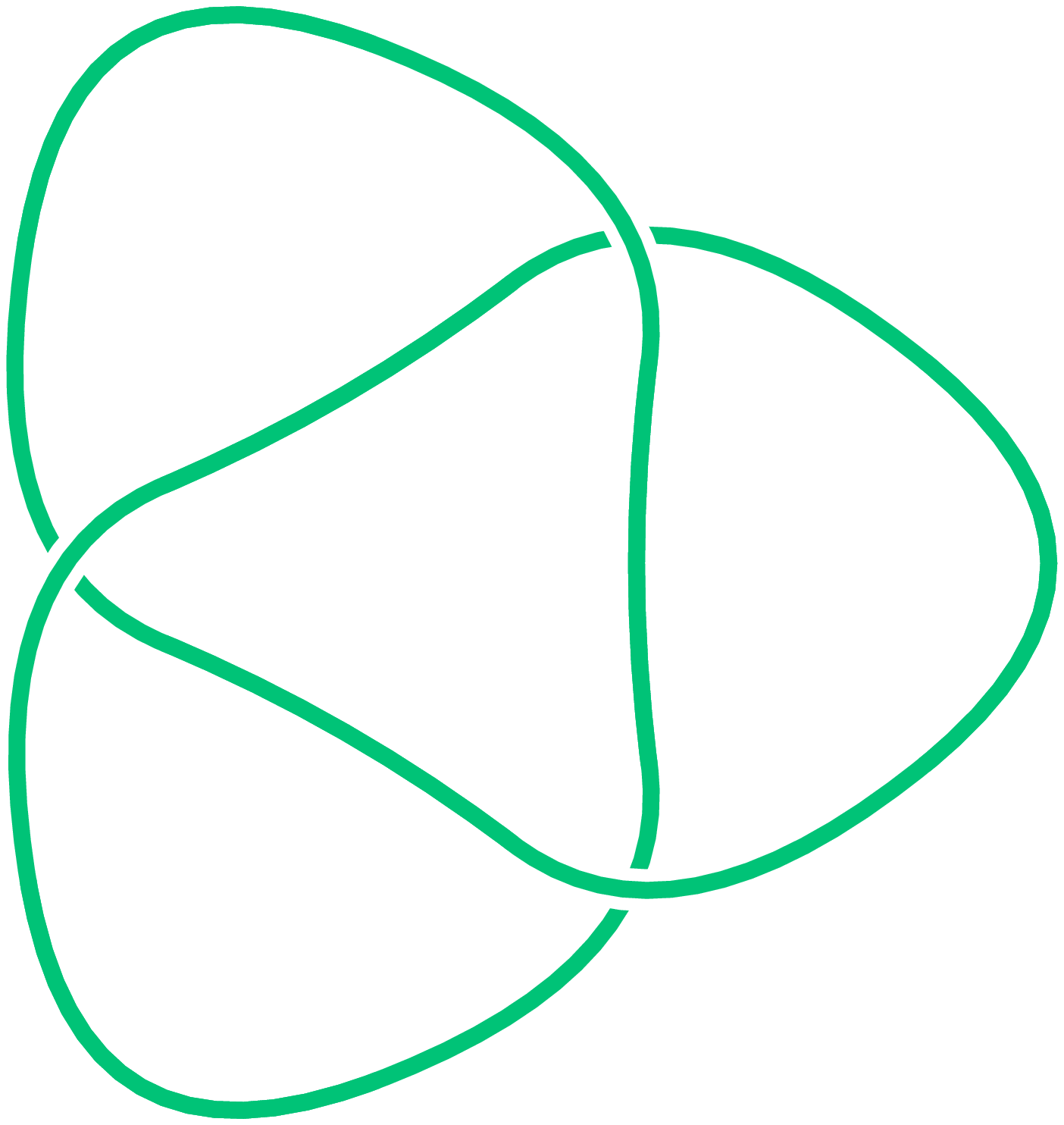}{1.3in} \qquad \psdraw{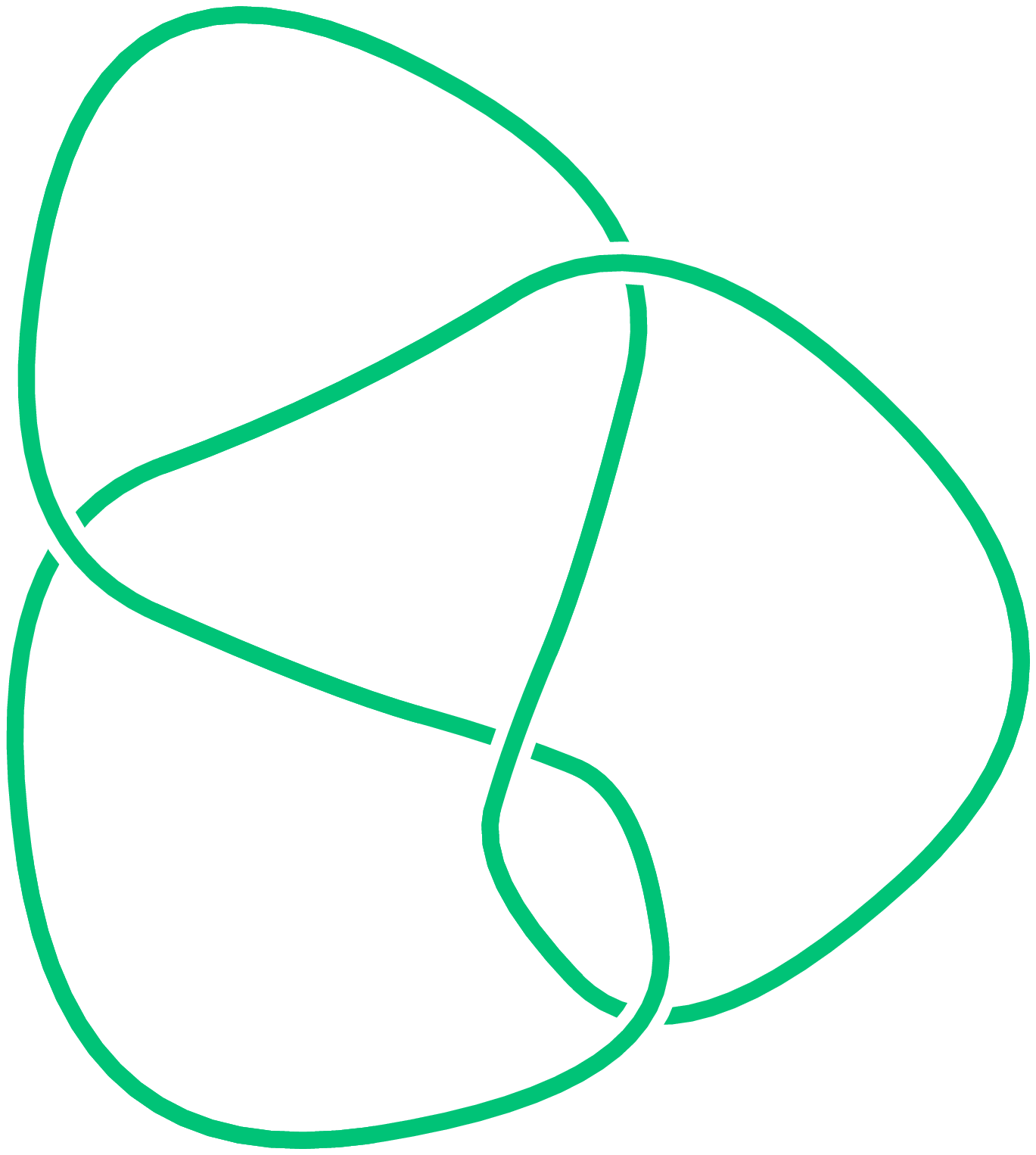}{1.3in}
$$

By the very definition, knots are flexible objects defined up to isotopy,
which allows the embedding to move in a smooth and arbitrary way as long as
it does not cross itself. In algebraic topology, a common way of studying 
knots (and more generally, spaces) is to associate computable numerical
invariants (such as Euler characteristic, or Homology). Invariants are
useful in deciding whether two knots are not the same. It is a much harder
problem to construct computable invariants that separate knots.

The invariant that we will consider in this paper is the Jones polynomial
of a knot; \cite{J}, which is a Laurent polynomial with integer coefficients,
associated to each knot. The quantum nature of the Jones polynomial is 
apparent both in the original definition of Jones
(using Temperley-Lieb algebras) and in the reformulation, due to Witten,
in terms of the expectation value of a Quantum Field Theory; see \cite{J,Wt}.

The combinatorics associated to a planar projection of a knot show that
the Jones polynomial is a computable invariant. However, it is hard to see
from this point of view the relation between the Jones polynomial and Geometry.
In Quantum Field Theory, one often reproduces Geometry by moving carefully
chosen parameters of the theory to an appropriate limit.

In our case, we will introduce a new parameter, a natural number which roughly
speaking corresponds to taking a connected $n$-parallel of a knot.
The resulting invariant is no longer a Laurent polynomial, but rather 
a sequence of Laurent polynomials.

The {\em colored Jones function} of a knot $K$ in 3-space 
is a sequence of Laurent polynomials 
$$
J_K: \BN \longto \BZ[q^{\pm}].
$$
The first term in the above sequence, $J_K(1)$ is the Jones polynomial
of $K$; see \cite{GL1}.

\subsection{The Hyperbolic Volume Conjecture}
\lbl{sub.hyperbolic}

Although knots are flexible objects, Thurston had the idea that their
complements have a unique decomposition in pieces of unique ``crystaline''
 shape.
The shapes in question are the $8$ different geometries in dimension
$3$, and the idea in question was termed the ``Geometrization Conjecture''.
The most common of the $8$ geometries is Hyperbolic Geometry, that is
the existence of a complete, finite volume, constant curvature $-1$ 
Riemannian metric on knot complements. Thurston proved that unless the
knot is torus or a satellite, then it carries a unique such metric; 
see \cite{Th}.

The {\em Hyperbolic Volume Conjecture} (HVC, in short) connects two very
different views of knot: namely Quantum Field Theory and Riemannian
Geometry. The HVC states for every hyperbolic knot $K$
$$
\lim_{n \to \infty} \frac{\log|J_K(n)(e^{\frac{2 \pi i}{n}})|}{n}=
\frac{1}{2 \pi}\, \text{vol}(S^3-K).
$$
where $\text{vol}(S^3-K)$ is the {\em volume} of a complete hyperbolic metric 
in the knot complement $S^3-K$. The conjecture was formulated in this form
by Murakami-Murakami \cite{MM} following an earlier version due to 
Kashaev, \cite{K}. More generally, Gukov (see \cite{Gu})
formulated a Generalized Hyperbolic Volume Conjecture that identifies the
limit
$$
\lim_{n \to \infty} \frac{\log|J_K(n)(e^{\frac{2 \pi i\a}{n}})|}{n}
$$
of a hyperbolic knot with known hyperbolic invariants (such as the volume
of cone manifolds obtained by hyperbolic Dehn filling), for $\a \in (0, 1]-\BQ$
or $\a=1$. Actually, the GHVC is stated for {\em complex} numbers $\a$.
For simplicity, we will study asympotics for real $\a \in [0,1]$.

At present, it is not known whether the limit in the HVC exists, let alone
that it can be computed. Explicit finite multisum formulas for the colored
Jones function of a knot exist; see for example \cite{GL1}. From these formulas
alone, it is difficult to study the above limit. In a sense, the question
is to understand the sequence of Laurent polynomials that appears in the
HVC. If the sequence is in some sense random, then it is hard to expect that
the limit exists, or that it can be computed.

Since the first term of this sequence is the Jones polynomial, and since
we know little about the possible values of the Jones polynomial, one would
expect that there is even less to be said about the colored Jones function.

\subsection{$q$-difference equations}
\lbl{sub.q}

Luckily, the colored Jones function behaves in a better way than its first
term, namely the Jones polynomial. This can be quantified by recent work
of TTQ Le and the first author, who proved that the colored Jones
function of a knot satisfies a $q$-difference equation. 

In other words, for every knot $K$ there exist rational functions $b_1(u,v),
\dots, b_d(u,v) \in \BQ(u,v)$ (which of course depend on $K$) such that
for all $n \in \BN$ we have:
$$
\sum_{j=0}^d b_j(q^n,q) J_K(n+j)=0.
$$

This opens the possibility of studying the $q$-difference equation rather
than one of its solutions, namely the colored Jones function.
Although the $q$-difference equation is not unique, it was shown by the
first author in \cite{Ga1} that one can choose a unique $q$-difference
equation, which is a knot invariant. Moreover, it was conjectured in
\cite{Ga1} that the characteristic polynomial of this
$q$-difference equation determines the characters of $\mathsf{SL}_2(\BC)$
representations of the knot complement, viewed from the boundary.

As was explained by the first author on several occasions, asymptotics
of solutions of $q$-difference equations would have consequences on the HVC.

In this introductory article we review the history of asymptotics of solutions
of $q$-difference equations.

\subsection{Asymptotics of differential
equations with a parameter}
\lbl{sub.bhistory}

Excellent references for differential equations with a parameter are
Olver's and Wasau's books; \cite{O} and \cite{Wa}. 
In 1837, Liouville and Green independently
studied systematically existence of formal (i.e., perturbative)
and actual solutions for second 
order differential equations with a parameter; see \cite{Gr,L}. 
Second order equations are very important for classical and quantum physics.

In 1908 Birkhoff had the insight to 
introduce and study {\em arbitrary} order differential equation with
a parameter (see \cite{B1}):
\begin{equation}
\lbl{eq.bir}
y^{(n)} + \rho a_{n-1}(x,\rho) y^{(n-1)}
+ \dots + \rho^n a_{0}(x,\rho) y=0
\end{equation}
where $y=y(x,\rho)$, 
$y^{(n)}$ means $n$-th derivative with respect to $x$ 
(assumed to be restricted to a real interval), and 
$\rho$ is a large complex parameter, and where the coefficient
$a_j(x,\rho)$ are complex $C^{\infty}$ functions with an expansion
$$
a_j(x,\rho)=\sum_{s=0}^\infty a_{j,s}(x) \rho^{-s}
$$
Birkhoff's working assumption was that the eigenvalues $\l_1(x),\dots,
\l_n(x)$ of the characteristic equation
$$
\l^n + a_{n-1,0}(x) \l^{n-1} + \dots + a_{0,0}(x) =0
$$
were distinct but not necessarily nowhere vanishing.

In 1926, three theoretical physicists,  
Wentzel-Krammer-Brillouin studied the second order
differential equation \eqref{eq.bir} under the assumption that its eigenvalues
do not collide, and developed connection formulas linking solutions in the 
exponential region with those in the oscillatory region. Their method is 
often referred to as the WKB method.

\subsection{Asymptotics of difference
equations}
\lbl{sub.bhistory2}

As a motivation for our
results, let us recall some fundamental results
of Birkhoff and Trjitzinsky from 1930 on difference equations without
a parameter; see \cite{B2}
and \cite{BT}.

A {\em difference equation} for a discrete function $f:\BN\to\BC$ has the form:
\begin{equation}
\lbl{eq.diff}
\sum_{j=0}^d a_j(n) f(n+j)=0
\end{equation}
where $a_j: \BN \longto \BC$ are discrete functions so that $a_0(n) a_d(n) 
\neq 0$ for all $n \in \BN$.
We will assume the existence of asymptotic expansions of $a_j(n)$ around
$n \to \infty$ for all $j=1,\dots,d$:
$$
a_j(n) \sim_{n \to \infty} n^{d_j/\w}(a_{j,0} + a_{j,1} n^{-1/\w} +
a_{j,2} n^{-2/\w} + \dots )
$$
where $\w \in \BN$. This certainly holds for $\w=1$
if $a_j$ are rational functions of $n$, as is often the case in combinatorial
problems.

Due to the nowhere vanishing of $a_d \cdot a_0$, it follows that the set  
of solutions of \eqref{eq.diff} is a vector space of dimension $d$.

There are two main problems of difference equations:
\begin{itemize}
\item
Existence of {\em formal series solutions} $\tpsi_1,\dots,\tpsi_d$
to \eqref{eq.diff}.
\item
Existence of a basis $\{\psi_1, \dots, \psi_d\}$ of solutions so that
$\psi_k(n)$ is asymptotic to $\tpsi_k(n)$ for large $n$.
\end{itemize}

In \cite{B2}, Birkhoff solved the existence of formal solutions in complete
generality (that is, without any assumptions on the eigenvalues
of the characteristic equation). 
In a sequel paper \cite{BT}, Birkhoff-Trjitzinsky solved the
second problem in complete generality.

Among other things, the formal solutions of Birkhoff lead to the development
of differential Galois theory, see \cite{vPS}.

Decades later, the results of Birkhoff and Trjitzinsky on difference
equations have found applications to enumerative combinatorics 
and numerical analysis; 
see for example Wimp and Zeilberger in \cite{Wi,WZ} and references therein.
It is not surprising that difference equations are used in numerical analysis,
since difference equations are numerical schemes
of approximating differential equations. 
In enumerative combinatorics and complexity theory, 
difference equations appear in recursive computation.
For example, the number $f(n)$ of involutions of $\{1,2,\dots,n\}$
(that is, permutations which are a product of 1 and 2-cycles) is given by
$$
f(n+2)=f(n+1) + (n+1) f(n)
$$
with $f(1)=1$, $f(2)=2$. Using the results of Birkhoff-Trjitzinsky and the
fact that $f(n)$ is monotone, it follows that
$$
f(n) \sim_{n\to\infty} K n^{n/2} e^{-n/2 + n^{1/2}} \left(1 + \frac{c_1}{
n^{1/2}} + \frac{c_2}{n} + \frac{c_3}{n^{3/2}} + \dots \right)
$$
for nonzero constants $c_i$ and some $K >0$. Actually, the $c_k$ can be 
computed recursively from the difference equation; see \cite[p.169]{WZ}.

\subsection{Asymptotics of difference equations with a parameter}
\lbl{sub.asymptotics}

By some historical accident, asymptotics of solutions of difference equations
with a parameter was not discussed a century ago. The first paper
that discusses second order difference equation with a parameter appears
to be due to Deift-McLaughlin (see \cite{DM}) which was generalized by
Costin-Costin to arbitrary order difference equations, \cite{CC}.

The  purpose in this paper is to show that for regular $q$-difference 
equations, a regular solution 
has a well-defined and computable exponential growth rate in terms of
a relative entropy of the characteristic polynomial of the 
$q$-difference equation; see Theorem \ref{thm.1} below.

This subject is classical and has been reinvented over the past hundred
years by several groups, often unaware of each others results.
In a sense, the problem of formal solutions of $q$-difference equations
is a problem in differential Galois theory; \cite{vPS}, 
and a problem in numerical analysis; see for example \cite{CC}. 

Our results are hardly new and are contained or can be obtained by
minor modifications from results of Costin-Costin or from work of Birkhoff and 
collaborators, \cite{B1,BT,CC}.

Since the presentation in the above papers varies by time and taste, we have
decided to give a self-contained account of the theory with complete proofs.
Hopefully, this will benefit the researchers in Quantum Topology and in
Analysis.

\subsection{Statement of the results}
\lbl{sub.qdiff}

In this paper, we will describe asymptotics of solutions of 
$q$-difference equations.

A {\em $q$-difference equation} for a sequence $(f(1),f(2),f(3),\dots)$ 
of smooth 
functions of $q$ has the form:
\begin{equation}
\lbl{eq.qdiff}
\sum_{j=0}^d b_j(q^k,q) f(k+j,q)=0
\end{equation}
where $b_j(v,u)$ are smooth functions  and $f(k,q)=f(k)(q)$.

Before we proceed further, let us remark that $q$ is a {\em variable} in 
\eqref{eq.qdiff} and not a complex number of absolute value less
(or more) than $1$. In the usual analytic theory of $q$-difference equations,
$q$ is a complex number inside or outside the unit disk.

Moreover, in the GHVC, we need to compute the $n$th term $f(n,q)$
in the above $q$-difference equation, and then evaluate it at 
$q_n=e^{2 \pi i \a/n}$, for $\a$ fixed. 
In other words, in the GHVC, $q_n$ is a complex number
that varies with $n$ in such a way that it stays in the unit circle and
approaches $1$ as $n \to \infty$.

With this in mind, $\e$-difference equations (defined below) are obtained
from $q$-difference equations by the substitution $q=e^{2 \pi i \e}$
where $\e$ is a small nonnegative real number, that plays the role
of Planck's constant.

The {\em characteristic polynomial} of the $q$-difference equation 
\eqref{eq.qdiff} is
$$
P(v,\l)=\sum_{j=0}^d b_j(v,1) \l^j
$$

\begin{definition}
\lbl{def.regularq}
We will say that \eqref{eq.qdiff} is {\em regular} if
$$
\mathrm{Dsc}_{\l} P(v,\l) \cdot b_0(v,1) \cdot b_d(v,1) \neq 0
$$
for all $v \in S^1$, where $\mathrm{Dsc}_{\l} P(v,\l)$ is the 
{\em discriminant} of $P(v,\l)$, which is a polynomial in the coefficients 
of $P(v,\l)$. 
\end{definition}

Let $\l_1(v), \dots, \l_d(v)$ denote the roots of the characteristic 
polynomial,  which we call the {\em eigenvalues} of \eqref{eq.qdiff}.
It turns out that \eqref{eq.qdiff} is regular iff 
the eigenvalues $\l_1(v), \dots, \l_d(v)$ never collide and never vanish, 
for every $v \in S^1$. Moreover, it follows by the implicit function theorem 
that the roots are smooth functions of $v \in S^1$.

Since we are interested in asymptotics of solutions of $q$-difference 
equations which, as we shall see, are governed by the magnitude of the
eigenvalues, we need to partition the circle  
according to the magnitudes
of the eigenvalues.

Let $S^1=\cup_{p \in \P} I_p$ denote a partition of $S^1$ into a finite union
of closed arcs (with nonoverlapping interiors), such that the magnitude
of the eigenvalues does not change in each arc. In other words,
for each $p \in \P$, there is a permutation $\s_p$ of the set $\{1,\dots,
d\}$ 
such that
$$
|\l_{\s_p(1)}(v)| \geq |\l_{\s_p(2)}(v)| \geq \dots \geq
|\l_{\s_p(d)}(v)| \qquad \text{for all} \qquad v \in I_p. 
$$

The following definition introduces a locally fundamental set of
solutions of $q$-difference equations.

\begin{definition}
\lbl{def.fundamentalsolq}

Fix a partition of $I$ as above. 
A set $\{\psi_1,\dots,\psi_d\}$ is a {\em locally fundamental set of solutions}
of \eqref{eq.qdiff} iff for every solution $\psi$ 
for every $p \in \P$ and for every $m=1,\dots,d$ 
there exist smooth functions $c^p_m$ such that
\begin{equation}
\lbl{eq.cpms}
\psi(k,q)=c^p_1(q)\psi_{\s_p(1)}(k,q) + \dots + c^p_d(q)
\psi_{\s_p(d)}(k,q)
\qquad \text{for all} \quad (k,q): \quad q^k \in I_p.
\end{equation}
\end{definition}

\begin{theorem}
\lbl{thm.1}
Assume that \eqref{eq.qdiff} is regular. Then, there exists a locally
fundamental set of solutions $\{\psi_1,\dots,\psi_d\}$ such that
\begin{itemize}
\item
For every $m=1,\dots,d$ and $(k,q)$ such that $q^k \in I_p$ we have
$$
\psi_m\left(k,e^{\frac{2 \pi i \a}{n}}\right)=
\exp\left(\frac{n}{\a} \Phi_m\left(\frac{k\a}{n},\frac{\a}{n}\right)\right).
$$
\item
for some smooth functions $\Phi_m$ with
uniform (with respect to $x \in I=[0,1]$) {\em asymptotic expansion}
$$
\Phi_m(x,\e) \sim_{\e \to 0}
\sum_{s=0}^\infty  \phi_{m,s}(x)\e^s 
$$
where $\phi_{m,s} \in C^\infty(I)$ for all $s$ 
\item
and leading term
\begin{equation}
\lbl{eq.int}
\phi_{m,0}(x)=\int_{0}^{x} \log(\l_m(e^{2 \pi i t})) dt
\end{equation}
where we have chosen a branch for the logarithm of $\l_m$.
\end{itemize}
\end{theorem}

\begin{remark}
\lbl{rem.determine}
For every $j=1,\dots,d$ the smooth functions $\phi_{j,s}$ for positive $s$
are uniquely determined from the coefficients $b_j(u,v)$ of \eqref{eq.qdiff}
by a hierarchy of first-order differential equations along with specified initial conditions. On the other hand,
the smooth functions $\Phi_m$ are not uniquely determined, since they are
obtained by a smooth interpolation. Thus, the locally fundamental set of 
solutions is not uniquely determined from the $q$-difference equation,
although its asymptotic behavior is.
\end{remark}

It follows from Theorem \ref{thm.1} that each locally fundamental solution 
$\psi_m(n,q)$ of the $q$-difference equation \eqref{eq.qdiff} satisfies the
GHVC in the sense that for every $\a \in [0,1]$ we have:  
$$
\lim_{n \to \infty} 
\frac{\log|\psi_m(n,e^{\frac{2 \pi i\a}{n}})|}{n}=
\int_{0}^{1} \log|\l_m(e^{2 \pi i \a t})| dt
$$

Fix a solution $\psi$ of \eqref{eq.qdiff}. Theorem \ref{thm.1} expresses
$\psi$ as a linear combination of $\psi_m$'s in each arc $I_p$. For
every $p \in P$, let
\begin{equation}
\lbl{eq.Sp}
S_p=\{m \in \{1,\dots,d\} | \,\,\, c^p_m \neq 0\}.
\end{equation}

Later (in Section \ref{sub.assregularq})
we will define the notion of a regular solution to a $q$-difference
equation. 

As a prototypical example, consider an $q$-difference equation that satisfies
$|\l_1(v)| > |\l_j(v)|$ for all $j \neq 1$ and all $v \in S^1$. Then, any 
solution that satisfies $c_1(0) \neq 0$ (or more generally, $c_1$ has
a nonvanishing derivative at $0$) is regular.

\begin{remark}
\lbl{rem.stokesq}
It is possible that $S_p \neq S_{p+1}$. In other words, the restriction of
a fixed solution $\psi$ to different intervals $I_p$ may be a linear
combination of different $\psi_j$s. 
This is an important phenomenon, referred by the
name of {\em Stokes phenomenon}; see \cite{Wa}. 
\end{remark}

Our next definition captures the growth rate of regular solutions to
regular $q$-difference equations.

\begin{definition}
\lbl{def.Sentropyq}
Fix a collection $S=\{S_p|\, p \in \P\}$ of subsets of $\{1,\dots, d\}$. 
The $S$-{\em entropy} 
$$
\s_S: [0,1] \to \BR
$$
of the $q$-difference equation \eqref{eq.qdiff} is defined by
$$
\s_S(\a)=
\int_0^{1} \log \chi_S(e^{2 \pi i \a t})  dt.
$$
where $\chi_S: [0,1] \to \BR$ is defined by 
$$
\chi_S(v)=\max_{j \in S_p} |\l_{\s_p(j)}(v)| \qquad \text{if} \qquad
v \in I_p.
$$
The {\em entropy} of \eqref{eq.qdiff} is the set of functions
$$
\{ \s_S:[0,1]\to\BR \quad | \quad S \subset \{1,\dots, d\} \}.
$$
\end{definition}

Notice that the entropy of a $q$-difference equation is not a real number,
but rather a finite collection of functions.

The main result is the following

\begin{theorem}
\lbl{thm.assq}
If $f$ is an $S$-regular solution of the regular $q$-difference equation
\eqref{eq.qdiff}, 
then for every $\a \in [0,1]$ we have:
$$
\lim_{n \to \infty} \frac{\log|f(n)(e^{\frac{2 \pi i\a}{n}})|}{n}=
\s_S(\a)
$$
\end{theorem}

Finally, let us define the $J$-entropy of a knot. In \cite{Ga1} the first
author showed that to every knot $K$ one can associate a canonical
$q$-difference equation of degree $d$, and a specific solution of it, 
namely the colored Jones function of $K$.

The $q$-difference equation itself is an invariant of a knot, which (by
definition) is determined by the colored Jones function of the knot.
Thus, any invariant of the $q$-difference equation is also an invariant
of a knot, which is determined by the colored Jones function of the knot.

\begin{definition}
\lbl{def.Jentropy}
The $J$-entropy of a knot is the entropy 
of its associated $q$-difference 
equation. We denote the $J$-entropy of a knot $K$ by
$$
\{ \s^J_{S,K}:[0,1]\to\BR \quad | \quad S \subset \{1,\dots, d\} \}.
$$
\end{definition}

\subsection{What's next?}
\lbl{sub.next}

The paper was written in the spring of 2004. Since then, a number of 
papers that discuss the asymptotics of the colored Jones function have
appeared; see \cite{Ga2, Ga3, GL2, GL3}.

\subsection{Acknowledgement}
The results of this paper were announced in the JAMI 2003 meeting in Johns
Hopkins. The first author wishes to thank J. Morava for the invitation,
and P. Deligne who suggested the asymptotic behavior of solutions of
$q$-difference equations.
The first author wishes to thank D. Boyd for sharing and explaining
his unpublished work and also A. Riese, T. Morley,  
and D. Zeilberger.

\section{$\e$-difference equations}
\lbl{sec.cd}

\subsection{$\e$-difference equations}
\lbl{sub.cd}

In this section, we will translate asymptotics of solutions of
$q$-difference equation in terms of asymptotics of solutions of
{\em $\e$-difference equations}. The latter are defined as follows.

Fix a positive number $\e_0$, a compact interval $I$ of $\BR$ and a natural
number $d$. 
We will consider functions $\phi: \D_{\e_0,I}\to \BC$ with domain 
\begin{equation}
\lbl{eq.domain}
\D_{\e_0,I}:=
\{(k\e,\e)| \quad k \in \BN, \e \in (0,\e_0], \quad
k \e \in I \}.
\end{equation}

Consider the $\e$-difference equation
for a function $\phi:\D_{\e_0,I}\to\BC$ 
\begin{equation}
\lbl{eq.ediff}
\sum_{j=0}^d a_j(k\e,\e) \phi((k+j)\e,\e)=0
\end{equation}
where $a_j \in C^\infty(I \times [0, \e_0])$.

We will assume that for all $j=0, \dots, d$, $a_j(x,\e)$ has a uniformly
(with respect to $x$) asymptotic expansion
\begin{equation}
\lbl{eq.ajass}
a_j(x,\e) \sim_{\e \to 0} \sum_{s=0}^\infty a_{j,s}(x) \e^s
\end{equation}
where $a_{j,s} \in C^\infty(I)$.

As we mentioned before,
$\e$-difference equations are obtained
from $q$-difference equations by the substitution $q=e^{2 \pi i \e}$
where $\e$ is a small nonnegative real number, that plays the role
of Planck's constant.

The {\em characteristic polynomial} of \eqref{eq.ediff} is
$$
P(x,\l)=\sum_{j=0}^d a_j(x,0) \l^j
$$

\begin{definition}
\lbl{def.regularc}
We will say that \eqref{eq.ediff} {\em regular} if
$$
\mathrm{Dsc}_{\l} P(x,\l) \cdot a_0(x,0) \cdot a_d(x,0) \neq 0
$$
for all $x \in I$.
\end{definition}

Let $\l_1(x), \dots, \l_d(x)$ denote the roots of the characteristic 
polynomial,  which we call the {\em eigenvalues} of \eqref{eq.ediff}.

It turns out that \eqref{eq.ediff} is regular iff 
the eigenvalues $\l_1(x), \dots, \l_d(x)$ never collide and never vanish, 
for every $x \in I$. Moreover, it follows by the implicit function theorem 
that the roots are smooth functions of $x \in I$.

Since we are interested in asymptotics of solutions of $\e$-difference 
equations which, as we shall see, are governed by the magnitude of the
eigenvalues, we need to partition  the interval $I$ according to the magnitudes
of the eigenvalues.

Let $I=\cup_{p \in \P} I_p$ denote a partition of $I$ into a finite union
of closed intervals (with nonoverlapping interiors), such that the magnitude
of the eigenvalues does not change in each interval. In other words,
for each $p \in \P$, there is a permutation $\s_p$ of the set $\{1,\dots,
d\}$ 
such that
$$
|\l_{\s_p(1)}(x)| \geq |\l_{\s_p(2)}(x)| \geq \dots \geq
|\l_{\s_p(d)}(x)| \qquad \text{for all} \qquad x \in I_p. 
$$

The following definition introduces a locally fundamental set of
solutions of $\e$-difference equations.

\begin{definition}
\lbl{def.fundamentalsol}
Fix a partition of $I$ as above.
A set $\{\psi_1,\dots,\psi_d\}$ is a {\em locally fundamental set of solutions}
of \eqref{eq.ediff} iff for every solution $\psi:\D_{\e,I}\to\BC$, 
for every $p \in \P$ and for every $m=1,\dots,d$ 
there exist smooth functions $c^p_m \in C^\infty [0,\e]$ 
such that
$$
\psi(k\e,\e)=c^p_1(\e)\psi_{\s_p(1)}(k\e,\e) + \dots + c^p_d(\e)
\psi_{\s_p(d)}(k\e,\e)
\qquad \text{for all} \qquad (k\e,\e) \in \D_{\e,I}.
$$
Here, the notation $c^p_m$ does not indicate the $p$th power of $c_m$.
\end{definition}

The next theorem summarizes the results of Costin-Costin.

\begin{theorem}
\lbl{thm.2}(\cite{CC})
Assume that \eqref{eq.ediff} is regular. Then, there exists a positive 
$\e' \leq \e_0$ and a locally fundamental set of solutions
$\{\psi_1,\dots,\psi_d\}$ of \eqref{eq.ediff} such that
\begin{itemize}
\item
For every $m=1,\dots,d$ and $(k\e,\e) \in \D_{e',I}$ we have
$$
\psi_m(k\e,\e)=\exp\left(\e^{-1} \Phi_m(k\e,\e)\right).
$$
\item
for some smooth functions $\Phi_m \in
C^\infty(I \times [0,\e'])$ with
uniform (with respect to $x \in I$) {\em asymptotic expansion}
\begin{equation}
\lbl{eq.Phim}
\Phi_m(x,\e) \sim_{\e \to 0}
\sum_{s=0}^\infty  \phi_{m,s}(x) \e^s 
\end{equation}
where $\phi_{m,s} \in C^\infty(I)$ for all $s$ 
\item
and leading term
\begin{equation}
\lbl{eq.intc}
\exp(\phi_{m,0}'(x)) = \l_m(x) .
\end{equation}
\end{itemize}
\end{theorem}


Fix a solution $\psi$ of \eqref{eq.ediff}. Theorem \ref{thm.2} expresses
$\psi$ as a linear combination of the $\psi_j$'s 
in each interval $I_p$. 
For
every $p \in P$, let
\begin{equation}
\lbl{eq.Spp}
S_p=\{m \in \{1,\dots,d\} | \,\,\, c^p_m \neq 0\}.
\end{equation}

Later (in Section \ref{sub.assregulare})
we will define the notion of a regular solution to an $\e$-difference
equation. 

As a prototypical example, consider an $\e$-difference equation that satisfies
$|\l_1(x)| > |\l_j(x)|$ for all $j \neq 1$ and all $x \in I=[a,b]$. Then, any 
solution that satisfies $c_1(a) \neq 0$ (or more generally, $c_1$ has
a nonvanishing derivative at $a$) is regular.

Our next definition captures the growth rate of regular solutions to
regular $\e$-difference equations.

\begin{definition}
\lbl{def.Sentropy}
Fix a collection $S=\{S_p|\, p \in \P\}$ of subsets of $\{1,\dots, d\}$. 
The $S$-{\em entropy} 
$$
\s_S: I \to \BR
$$
of the $\e$-difference equation \eqref{eq.ediff} is defined by
$$
\s_S(x)=
\int_0^{x} \log \chi_S(t)  dt.
$$
where $\chi_S: I \to \BR$ is defined by 
$$
\chi_S(x)=\max_{j \in S_p} |\l_{\s_p(j)}(x)| \qquad \text{if} \qquad
x \in I_p.
$$
\end{definition}

\begin{theorem}
\lbl{thm.ass}
If $\psi$ is an $S$-regular solution to a regular $\e$-difference equation, 
and $x \in I$, we have:
$$
\lim_{\e\to 0^+} \e \log |\psi(x,\e)|=
\s_S(x)
$$
\end{theorem}

The next remarks concern the uniqueness of a set of locally fundamental
solutions to \eqref{eq.ediff}.

\begin{remark}
\lbl{rem.determine2}
For every $m=1,\dots,d$ the smooth functions $\phi_{m,s}$ for positive $s$
are uniquely determined by \eqref{eq.qdiff} and the initial condition
$\phi_{m,s}(0)=0$.
Indeed, applying Taylor series (with respect to $\e$) in \eqref{eq.ediff}
and collecting terms, we get for example: 
\begin{eqnarray*}
\phi_{m,1}'(x) &=& 
-\frac{1/2 \phi_{m,0}''(x)  \sum_{j=0}^d a_j(x,0) j^2 \l_m^j(x)
+ \sum_{j=0}^d \pt_{\e} a_j(x,0) \l_m^j(x)}{
\sum_{j=0}^d a_j(x,0) j \l_m^j(x)} \\
&=& -\left(\frac{1}{2} \l' \frac{P_{\l\l}}{P_{\l}} +
\frac{P_{\e}}{P_{\l} \l} \right) \left|_{\l=\l_m(x)} \right.
\end{eqnarray*}
where $f_x(x,\l)=\pt/\pt_x f(x,\l)$ and $f_{\l}(x,\l)=\pt/\pt_{\l} f(x,\l)$.

Similarly, for $s \geq 1$ we have:
$$
\phi_{m,s}'(x)=-\frac{H_s(x)}{\sum_{j=0}^d a_j(x,0) j \l_m^j(x)}
$$
where $H_s(x)$ is a function of derivatives of $a_j(x,0)$ and $\phi_{m,t}$
for $t < s$.
Notice that the denominator vanishes nowhere since the roots do not collide
and do not vanish for every $x \in I$.
\end{remark}

\begin{remark}
\lbl{rem.analytic}
If the coefficients $a_j(x,\e)$ of the regular $\e$-difference equation 
\eqref{eq.ediff} are analytic functions, then the functions $\phi_{m,s}$
of Theorem \ref{thm.2} are also analytic, for every $m$ and $s$. This
follows by induction from the differential hierarchy which these functions
satisfy, and from the fact that the eigenvalues are analytic functions.
Even though $\phi_{m,s}$ is analytic for every $m$ and $s$, 
the series
$$
\sum_{s=0}^\infty \phi_{m,s}(x) \e^s
$$
is in general divergent, and the functions $\Phi_{m,s}$ of Theorem 
\ref{thm.2} are not analytic.
\end{remark}

\begin{remark}
\lbl{rem.smoothunique}
Even though the functions $\phi_{m,s}$ are uniquely determined by the
$\e$-difference equation, the smooth functions $F_m$ (and thus the
locally fundamental set of solutions $\psi_m$) are not uniquely determined
by the $\e$-difference equation. The problem is that smooth
interpolation is not unique. Recently developed ideas of exponentially
small corrections might construct a unique set of locally fundamental
solutions when the coefficients of \eqref{eq.ediff} are analytic functions.
We will elaborate on this in a separate occasion.
\end{remark}

\subsection{Converting $q$-difference equations to $\e$-difference
equations}
\lbl{sub.3implies2}

The translation of $q$-difference equations to $\e$-difference equations
is as follows. If $f$ satisfies the $q$-difference equation 
$$
\sum_{j=0}^d b_j(q^k,q) f(k+j,q)=0
$$
then set
$$
q=e^{2 \pi i \e}, \quad
b_j(e^{2 \pi i x} ,e^{2 \pi i \e})=a_j(x, \e)
$$
and consider the $\e$-difference equation for a function $\phi$ 
(with domain $\D_{\e_0,I}$ for some $\e_0>0$ and $I=[0,2\pi]$):
$$
\sum_{j=0}^d a_j(k\e,\e) \phi((k+j )\e,\e)=0
$$

The following lemma, although elementary, is the key to translating
$q$-difference equations to $\e$-difference equations.

\begin{lemma}
\lbl{lem.translate}
For every $(k\e,\e) \in I \times [0,\e_0]$ we have:
\begin{equation}
\lbl{eq.phif}
\phi(k\e,\e)=f(k,e^{2 \pi i \e}).
\end{equation}
Consequently, for every $\a \in (0,1]$, we have:
\begin{equation}
\lbl{eq.phiff}
\lim_{k\to \infty} \frac{1}{k} \log|f(k,e^{2 \pi i\a/k})| =\a^{-1}
\lim_{\e\to 0}\e \log |\phi(\a,\e)|.
\end{equation}
Thus, Theorem \ref{thm.2} implies Theorem \ref{thm.1}.
\end{lemma}

\begin{proof}
Observe that $a_j(k\e,\e)=b_j(e^{2 \pi i k\e}, e^{2 \pi i \e})$, thus
$(k,\e)\to\phi(k\e,\e)$ satisfies the equation
$$ 
\sum_{j=0}^d b_j(e^{2 \pi i k\e}, e^{2 \pi i \e}) \phi((k+j)\e,\e)=0
$$
and so does $(k,\e)\to f(k, e^{2 \pi i \e})$. Since solutions with matching
initial conditions are unique, \eqref{eq.phif} follows.

Equation~\eqref{eq.phiff} follows from equation~\eqref{eq.phif} by the substitution $\e=\a/k$:
$$
\frac{1}{k} \log  |f(k,e^{2 \pi i\a/k})| =
\frac{1}{k} \log  |\phi(k\e,\a/k)|=
\a^{-1} \e \log |\phi(\a,\e)|.
$$
\end{proof}

\section{Some linear algebra}
\lbl{sec.linear}

In this section we will review some linear algebra.
It is obvious that the complex roots of a monic polynomial uniquely
determine it. It is also known \cite{GLR} that the eigenvalues of a companion matrix
uniquely determine it, in case they are distinct.

Consider a {\em companion} $d$ by $d$ matrix 
$$
A=\begin{pmatrix}
0     & 1     & 0 & \dots & 0 \\
0   & 0   & 1   & \dots & 0   \\
0 & 0     & 0  & \dots & 0   \\
0 & 0     & 0 & \dots & \dots \\
\dots & \dots & \dots & \dots & 1 \\
-a_{0} & -a_{1} & -a_{2} & \dots & -a_{d-1} 
\end{pmatrix}
$$
The characteristic polynomial of $A$ is
$$
\l^d + \sum_{j=0}^{d-1} a_j \l^j
$$
with roots $\l_1,\dots,\l_d$. Let $M=(\l_j^{i-1})_{i,j}$ be the
Vandermonde matrix, and $D=\text{diag}(\l_1, \dots, \l_d)$
be the diagonal matrix with diagonal entries $\l_1, \dots, \l_d$.

\begin{lemma}
\lbl{lem.0}
If a companion matrix has distinct eigenvalues, then with the above notation
we have: 
$$
A=M D M^{-1}
$$ 
\end{lemma}

\begin{proof}
Observe that $v_j=(1,\l_j,\dots, \l_j^{d-1})^T$ is an eigenvector of $A$
with eigenvalue $\l_j$. Thus, $M=(v_1, \dots, v_d)$ and
$AM=MD$. The result follows.
\end{proof}

Now, consider a companion matrix $A(u)$ whose entries in the bottom row
are smooth functions in a variable $u$, 
with roots $\l_1(u),\dots,\l_d(u)$ which we assume are distinct for all $u$.

The next lemma is a key estimate for the norm of long products of slowly
varying matrices. In the language of physics, $A(u)$ is the {\em 
transfer matrix} and $A(n) \dots A(2)A(1)$ is the {\em transition matrix}.

\begin{lemma}
\lbl{lem.1}
Assume that the eigenvalues $\l_1(u),\dots,\l_d(u)$ of $A(u)$
are distinct for all $u$
and 
$$
\max_j \sup_u | \l_j(u)| \leq 1 + C \e.
$$ 
Then for $m \leq n$, $n \e \in I$, we have
$$
\Vert A(n) A(n-1) \dots A(m)\Vert \leq C'
$$
\end{lemma}

\begin{proof}
By Lemma \ref{lem.0}, we have:
$$
A(u)=M(u) D(u) M(u)^{-1}
$$ 
If $m \leq n$, it follows that 
$$
A(n) A(n-1) \dots A(m) =
M(n) D(n) M(n)^{-1} M(n-1) D(n-1) M(n-1)^{-1} \dots M(m) D(m) M(m)^{-1},
$$
which implies that
\begin{eqnarray*}
\Vert A(n) A(n-1) \dots A(m) \Vert & \leq &  
\Vert M(n) \Vert  M(m)^{-1} \Vert \cdot  
\Vert D(n) \Vert \dots \Vert D(m) \Vert \\
& & 
\Vert M(n)^{-1}M(n-1) \Vert \dots \Vert M(m+1)^{-1} M(m) \Vert.
\end{eqnarray*}
Now, we have 
\begin{eqnarray*}
\Vert D(k)\Vert & \leq & 1+ C\e \qquad \text{for} \quad k=m, \dots, n \\
\Vert M(k) M(k-1)^{-1}\Vert & \leq & 1 + C' \e \qquad \text{by Lemma 
\ref{lem.2}}.
\end{eqnarray*}
If $I=[a,b]$, using the fact that $n\e,m\e \in I$, we obtain: 
\begin{eqnarray*}
\Vert A(n) A(n-1) \dots A(m)\Vert & \leq & (1+C\e)^{n-m}(1+C'\e)^{n-m} \\
& \leq & (1 + C'' \e)^{2(n-m)} \\
& \leq & (1+C''\e)^{2\frac{b-a}{\e}} \\
& \leq & e^{2 C'' (b-a)}.
\end{eqnarray*}
\end{proof}

\begin{lemma}
\lbl{lem.2}
If $M=(x_j^{i-1})_{i,j}$ and $N=(y_j^{i-1})_{i,j}$ are Vandermonde
matrices, such that $M$ is nonsingular, then
$$
(M^{-1}N)_{i,j}= \prod_{l \neq i} \frac{y_j-x_l}{x_i-x_l}.
$$
\end{lemma}

\section{Existence of formal solutions}
\lbl{sub.formal}

In this section we will prove that \eqref{eq.ediff} has a unique set of
formal solutions. Let us define those first.

\begin{definition}
\lbl{def.formalseries}
A {\em formal series} $\tpsi(x,\e)$ is one of the form
\begin{equation}
\lbl{eq.eqf}
\tpsi(x,\e) =
\exp\left( \e^{-1} \sum_{s=0}^\infty  \phi_{s}(x) \e^s \right)
\end{equation}
where $\phi_{s} \in C^\infty(I)$ are smooth functions
for all $s$. 
\end{definition}

Note that $\e \log \tpsi(x,\e)$ lies in the ring $C^\infty(I)[[\e]]$ of
formal power series with coefficients smooth functions on $I$.

Note further that if $\tpsi(x,\e)$ is a formal series, so is
$\tpsi(x+j\e,\e)$ for every $j \in \BZ$, 
where the latter may be defined using the Taylor series
of $\phi_s(x+j\e)=\sum_{t=0}^\infty \frac{1}{t!} \phi^{(t)}_s(x) j^t \e^t$:
\begin{eqnarray*}
\lbl{eq.galois}
\tpsi(x+j\e,\e) &=&
\exp\left( \e^{-1}  \sum_{s=0}^\infty \left( \sum_{t=0}^s 
\frac{1}{t!} \phi_{s-t}^{(t)}(x) \right) \e^s  \right) \\ 
&=& \tpsi(x,\e)
\exp\left( \phi'_0(x)j + \left(\phi'_1(x)j + \frac{\phi''_0(x)}{2!}j^2\right)\e
+ \left(\phi'_2(x)j + \frac{\phi''_1(x)}{2!}j^2
+ \frac{\phi'''_0(x)}{3!}j^3 \right)\e^2 + \dots \right) 
\end{eqnarray*}

It follows that if $\tpsi(x,\e)$ is a formal series, then the ratio
$\tpsi(x+\e,\e)/\tpsi(x,\e)$ lies in the ring $C^\infty(I)[[\e]]$.

%

Using the language of {\em difference Galois theory} (see \cite[p.4]{vPS}) this
implies that 

\begin{lemma}
\lbl{lem.galois}
$C^\infty(I)[[\e]]$ is a finite difference ring, under the map $x\to x+\e$.
\end{lemma}

\begin{definition}
\lbl{def.formalsol}
We say that a formal series $\tpsi$ of \eqref{eq.eqf} is a 
{\em formal series solution} to \eqref{eq.ediff} iff 
\begin{equation}
\lbl{eq.fss}
\frac{1}{\tpsi(x,\e)}
\sum_{j=0}^d a_j(x,\e) \tpsi(x+j\e,\e)=0 \in C^\infty(I)[[\e]].
\end{equation}
\end{definition}
It is easy to see that if $\tpsi$ is a formal solution to \eqref{eq.ediff},
then the leading term $\phi_{0}$ satisfies the equation
\begin{equation}
\lbl{eq.intf}
\exp(\phi_{0}'(x)) = \l(x) 
\end{equation}
where $\l(x)$ is an eigenvalue of \eqref{eq.ediff}.

\begin{proposition}
\lbl{prop.formal}
Assume that \eqref{eq.ediff} is regular. Then, \eqref{eq.ediff}
has $d$ unique formal series solutions $\tpsi_1, \dots, \tpsi_d$ with
leading terms corresponding to the eigenvalues of \eqref{eq.ediff}.
\end{proposition}

\begin{proof}
First we need to show that \eqref{eq.fss} is indeed an equation in the
power series ring $C^\infty(I)[[\e]]$, i.e., that the terms involving negative
powers of $\e$ cancel.

Suppose that $\tpsi$ is given by \eqref{eq.eqf}. It follows from
the calculation preceding Lemma \ref{lem.galois} that 
for every $s \in \BN$, we have:
\begin{equation}
\lbl{eq.psie}
\text{coeff}\left(\e^s, \frac{\tpsi(x+j\e,\e)}{\tpsi(x,\e)}\right)=
\begin{cases}
\exp(j \phi_0'(x)) & \text{if} \quad s=0 \\
j \exp(j \phi_0'(x)) \phi'_s(x) + \text{terms}_s & \text{if} \quad s >0
\end{cases}
\end{equation}
where $\text{coeff}(\e^s,g(\e))$ denotes the coefficient of $\e^s$ in a power
series $g(\e)$, and 
where $\text{terms}_s$ is a polynomial in the 
derivatives of $\phi_t$ for $t < s$.

Expand the terms of Equation \eqref{eq.fss}
into power series in $\e$ using the above equation and \eqref{eq.ajass},
and collect terms of powers of $\e$. It follows that \eqref{eq.fss} is
equivalent to a hierarchy of first order differential equations:
\begin{eqnarray*}
\sum_{j=0}^d a_j(x,0) \exp(j \phi_0'(x)) &=& 0 \\
\sum_{j=0}^d a_j(x,0) j \exp(j \phi_0'(x)) \phi_s'(x) + \text{Terms}_s &=& 0
\end{eqnarray*}
where for positive $s$, $\text{Terms}_s$ is a   
polynomial in the derivatives of $\phi_t$ and $a_j(x,0)$ for $t < s$.

Now fix an $m \in \{1,\dots,d\}$, and choose $\phi_{m,0}$ such that
$\exp(\phi_{m,0}'(x))=\l_m(x)$, where $\l_1(x),\dots,\l_d(x)$ are the 
eigenvalues of \eqref{eq.ediff}.
Since \eqref{eq.ediff} is regular, it follows that the roots
$\l_1(x), \dots, \l_d(x)$ of the characteristic polynomial $P(x,\l)$
never collide, and never vanish for $x \in I$. Thus, 
$$
0 \neq \left(\l \frac{d}{d\l} P(x,\l)\right)_{\l=\l_m(x)} =
\sum_{j=0}^d j a_j(x,0) \l_m^j(x)
$$
for all $x \in I$. Thus, after we choose $\phi_{m,0}$, it follows that
we can find functions $\phi_{m,s}$ for $s \geq 0$ that satisfy the above
hierarchy. Moreover, for every $m$, the sequence of functions $\phi_{m,s}$
is uniquely determined by the above hierarchy and the initial conditions
$\phi_{m,s}(0)=0$.
\end{proof}

\subsection{An alternative formal series}
\lbl{sub.altformal}

In this section we present an alternative, and slightly more general form,
of formal series. In case of regular $\e$-difference equations this
alternative form will not be needed. However, when eigenvalues collide
or vanish, one must use this alternative form of formal series.
Thus, in the present paper we will not use this alternative form
of formal series, and the reader may skip this section.

\begin{definition}
\lbl{def.formalseriesalt}
An {\em alternative formal series} $\tpsi(x,\e)$ is one of the form
\begin{equation}
\lbl{eq.eqfalt}
\tpsi(x,\e) =
\exp\left( \e^{-1} \phi(x)\right) \sum_{s=0}^\infty  \phi_{s}(x) \e^s 
\end{equation}
where $\phi,\phi_{s} \in C^\infty(I)$ are smooth functions
for all $s$, and $\phi_0(x) \neq 0$ for all $x \in I$.
\end{definition}

In the remainder of this subsection, we will refer to alternative
formal series simply by formal series.

Note that if $\tpsi(x,\e)$ is a formal series, so is
$\tpsi(x+j\e,\e)$ for any $j \in \BZ$, 
where the latter may be defined using the Taylor series
of $\phi_s(x+j\e)$ and $\phi(x+j\e)$ around $x$. It follows that 

\begin{eqnarray*}
\lbl{eq.galoisalt}
\tpsi(x+j\e,\e) &=&
\exp\left( \e^{-1} \phi(x)\right)(\phi_0(x)+ ( \phi^{(1)}(x) \phi_0(x)
+ \phi^{(1)}_0(x) + \phi_1(x)) \e + \dots ). 
\end{eqnarray*}

Moreover, if $\tpsi(x,e)$ is a formal series, then the ratio
$\tpsi(x+\e,\e)/\tpsi(x,\e)$ lies in the ring $C^\infty(I)[[\e]]$.
This follows from \eqref{eq.galoisalt} and the following computation, valid
for every $j \in \BN$:

\begin{eqnarray*}
\tpsi(x+j\e,\e) &=& \exp(\e^{-1} \phi(x+j\e)) ( \phi_0(x+j\e)+ \phi_1(x+j\e)\e
+O(\e^2)) \\
&=&
\exp\left(\e^{-1} \phi(x) + \phi^{(1)}(x)j + \frac{\phi^{(2)}(x)}{2} j^2 \e +
O(\e^2)\right)
\left(\phi_0(x) + (\phi^{(1)}_0(x)j + \phi_1(x)) \e + O(\e^2)\right) \\
& = &
\tpsi(x,\e) \exp(\phi^{(1)}(x) j)\left(1 + 
\left(\frac{\phi^{(2)}(x)}{2} j^2 + \frac{\phi^{(1)}_0(x)}{\phi_0(x)}
\right)
\e + O(\e^2) \right)
\end{eqnarray*}

In analogy with Lemma \ref{lem.galois}, this implies that

\begin{lemma}
\lbl{lem.galoisalt}
$C^\infty(I)[[\e]]$ is a finite difference ring, under the map $x\to x+\e$.
\end{lemma}

\begin{definition}
\lbl{def.formalsolalt}
We say that a formal series $\tpsi$ of \eqref{eq.eqfalt} is a 
{\em formal series solution} to \eqref{eq.ediff} iff 
\begin{equation}
\lbl{eq.fssalt}
\frac{1}{\tpsi(x,\e)}
\sum_{j=0}^d a_j(x,\e) \tpsi(x+j\e,\e)=0 \in C^\infty(I)[[\e]].
\end{equation}
\end{definition}
It is easy to see that if $\tpsi$ is a formal solution to \eqref{eq.ediff},
then the leading term $\phi$ satisfies the equation
\begin{equation}
\lbl{eq.intfalt}
\exp(\phi'(x)) = \l(x) 
\end{equation}
where $\l(x)$ is an eigenvalue of \eqref{eq.ediff}.

In analogy with Proposition \ref{prop.formal}, we have the following:

\begin{proposition}
\lbl{prop.formalalt}
Assume that \eqref{eq.ediff} is regular. Then, \eqref{eq.ediff}
has $d$ unique formal series solutions $\tpsi_1, \dots, \tpsi_d$ with
leading terms corresponding to the eigenvalues of \eqref{eq.ediff}.
\end{proposition}

\section{Proof of Theorem \ref{thm.2}}
\lbl{sec.proof2}

In this section we prove Theorem \ref{thm.2}.
The strategy is to 
\begin{itemize}
\item[(a)]
prove that there exists a solution $\psi_1$ with the stated properties
where $\l_1(x)$ is an eigenvalue with maximum magnitude.
\item[(b)]
use this solution $\psi_1$ to reduce Theorem \ref{thm.2}
to the case of a $\e$-difference 
equation of degree one less than the original one.
\item[(c)]
prove that the constructed set of solutions is a locally fundamental set.
\end{itemize}

Without loss of generality, we will assume that the eigenvalues of
\eqref{eq.ediff} satisfy the inequality:
$$
|\l_1(x)|\geq |\l_2(x)| \geq \dots \geq |\l_d(x)|
$$
for all $x \in I$. Otherwise, we can partition $I$ into subintervals where 
this is true.

\subsection{Existence of a solution corresponding to the eigenvalue
of maximum magnitude}
\lbl{sub.exist}

Consider first a formal solution $\tpsi_1$ of \eqref{eq.ediff} given
in Proposition \ref{prop.formal}, which satisfies \eqref{eq.eqf} and
\eqref{eq.intf}. Consider the smooth functions $\phi_{1,s} \in C^\infty(I)$
of \eqref{eq.intf}.

The proof of the following lemma (due to Borel in case $\phi_{1,s}$ are
constant functions, for all $s$) can be found in \cite[Lemma 2.5]{GG}:

\begin{lemma}
\lbl{lem.borel}
There exists a smooth function $\hat\Phi_1 \in C^\infty(I \times [0,\e_0])$
 such that we have (uniformly in $x \in I$): 
$$
\hat\Phi_1(x,\e) \sim_{\e \to 0} \sum_{j=0}^\infty \phi_{1,s}(x) \e^s.
$$
\end{lemma}

Now, consider the unique solution $\psi_1$ of \eqref{eq.ediff} with initial
conditions 
$$
\psi_{1}(k\e,\e)=\exp(\e^{-1} \hat\Phi_1(k \e,\e))
\quad \text{ for } \quad k=0,\dots,d-1
$$
and for small enough $\e$, where without loss of generality, we assume
that $I=[0,b]$.

Of course, for large $k$ it may not be true that 
$\psi_{1}(k\e,\e)=\exp(\e^{-1} \hat\Phi_1(k \e,\e))$. The next proposition
estimates the error, uniformly with respect to $k$:

\begin{proposition}
\lbl{prop.estimate}
There exists an $\e'>0 $ and constants $C_s$ such that for all 
$(k\e,\e) \in \D_{\e',I}$ and all $s \in \BN$, we have (uniformly in $k$):
\begin{equation}
\lbl{eq.estimate}
\left| \frac{\psi_1(k\e,\e)}{\exp(\e^{-1} \hat\Phi_1(k \e,\e))}
-1 \right| < C_s \e^s
\end{equation}
\end{proposition}

\begin{proof}
Let us make a change of variables:
\begin{equation}
\lbl{eq.change}
\th=\frac{\psi_1}{\hpsi_1},
\end{equation}
where
$$
\hpsi_1(x,\e)=\exp\left(\e^{-1} \hat\Phi_1(x,\e) \right).
$$
We will show that for a fixed $s_0$, and for every $s$ there exists a constant
$C_s$ such that for all $(k\e,\e) \in \D_{\e_0,I}$ we have:
\begin{equation}
\lbl{eq.step1}
\left| \th(k\e,\e) -1 \right| < C_s \e^{s+1-s_0},
\end{equation}

Since  $\psi_1$ satisfies \eqref{eq.ediff}, 
it follows that $\th$ satisfies
\begin{equation}
\lbl{eq.ediff2}
\sum_{j=0}^d b_j(k\e,\e) \th((k+j)\e,\e)=0
\end{equation}
where 
$$
b_j(x,\e)=a_j(x,\e) \frac{\hpsi_1(x+j\e,\e)}{\hpsi_1(x,\e)}.
$$ 
It is easy to see that 
\begin{itemize}
\item
$b_j(x,\e) \in C^\infty(I \times [0,\e_0])$, has uniform (with respect
to $x$) $\e$-asymptotic expansion as in \eqref{eq.ajass},
\item
$b_{j,s}(x,0)=a_j(x,0) \l_1^j(x)$,
\item
and since $\tpsi$ is a formal series solution to \eqref{eq.ediff} and
$\hat\Phi_1$ is given by Lemma \ref{lem.borel}, if follows that for every $s$
we have:
\begin{equation}
\lbl{eq.approximate}
\sum_{j=0}^d b_j(x,\e)=O(\e^s)
\end{equation}
\end{itemize}

The characteristic polynomial of \eqref{eq.ediff2}
has roots $\mu_m(x):=\l_m(x)/\l_1(x)$ for $j=1,\dots,d$. 
If \eqref{eq.ediff} is regular, so is \eqref{eq.ediff2}.
Notice that $\l_1(x)$ vanishes nowhere since \eqref{eq.ediff2} is 
regular.

We now show \eqref{eq.step1}. 
Let us write the difference equation \eqref{eq.ediff2} in matrix form
\begin{equation}
\lbl{eq.matrixform}
\Th((k+1)\e,\e)=A(k\e,\e) \Th(k\e,\e)
\end{equation}
where
$$
\Th(x,\e)=
\begin{pmatrix}
\th(x,\e) \\ \th(x+\e,\e) \\ \th(x+2\e,\e) \\ \dots \\ \dots \\ 
\th(x+(d-1)\e,\e)
\end{pmatrix}
\qquad
A(x,\e)=\begin{pmatrix}
0     & 1     & 0 & \dots & 0 \\
0   & 0   & 1   & \dots & 0   \\
0 & 0     & 0  & \dots & 0   \\
0 & 0     & 0 & \dots & \dots \\
\dots & \dots & \dots & \dots & 1 \\
-c_{0}(x,\e) & -c_{1}(x,\e) & -c_{2}(x,\e) & \dots & -c_{d-1}(x,\e) 
\end{pmatrix}
$$
and $c_j(x,\e)=b_j(x,\e)/b_d(x,\e)$.
Iterating, we obtain that 
$$
\Th(k\e,\e)=A(k\e,\e) A((k-1)\e,\e) \dots A( \e,\e) \Th(0,\e)
$$
for all $k \geq 1$, where $\Th(0,\e)=\mathbb{1}$, 
a column vector with all entries equal to $1$.

Equation \eqref{eq.approximate} gives:
\begin{equation}
\lbl{eq.gives}
A(k \e,\e) \mathbb{1}= \mathbb{1} +  \e^s \mathbb{E}_{k}(\e)  
\end{equation}
where $\Vert\mathbb{E}_{k}(\e) \Vert < C$ uniformly in $k$ and $\e$.
Feeding in the above equation, we obtain:
\begin{equation}
\lbl{eq.th1}
\Th(k\e,\e)=\mathbb{1} + \e^s \mathbb{E}_k(\e) + \e^s
\sum_{j=1}^{k-1}  A(k\e,\e) A((k-1)\e,\e)
\dots A((j+1)\e,\e) \mathbb{E}_j(\e). 
\end{equation}

Now, let us look at the roots $\mu_1(x,\e), \dots, \mu_d(x,\e)$ of
$$
\sum_{j=0}^d b_j(x,\e) \mu^j=0.
$$
Since $\max_{j}\sup_{x \in I} |\mu_j(x,0)|=1$, it follows that
$$
\max_j \sup_{x \in I} |\mu_j(x,\e)| \leq 1 + C \e.
$$
Since $k\e$ lies in $I$, a compact set,
Lemma \ref{lem.1} and Equation \eqref{eq.th1} imply that 
\begin{equation}
\lbl{eq.imply}
\Vert\Th_k(\e) -\mathbb{1}\Vert \leq k C'_s \e^{s+1} \leq C_s \e^{s-s_0+1}
\end{equation}
for all $k\e \in I$, 
where $s_0=1$. This completes the proof of \eqref{eq.step1}.

Equations \eqref{eq.step1} and \eqref{eq.estimate} differ in the presence
of $s_0$. It is easy to see that if $f$ is a function such that for a fixed
$s_0$ and any $s \in \BN$ we have:
$$
| f(\e)-\sum_{t=0}^{s+s_0} C_t \e^t | < D_s \e^{s+1},
$$
then
$$
| f(\e)-\sum_{t=0}^{s} C_t \e^t | < (D_s + |C_{s+1}|+\dots |C_{s+s_0}|) 
\e^{s+1}.
$$
This observation shows that \eqref{eq.step1} implies \eqref{eq.estimate}
and concludes the proof of the Proposition.
\end{proof}

\begin{proposition}
\lbl{prop.step1}
$\mathrm{(a)}$ There exists a smooth function $\Phi_1 \in C^\infty(I \times 
[0,\e'])$ such that
\begin{itemize}
\item[(a)] For all $(k\e,\e) \in \D_{\e',I}$, we have:
$$
\psi_1(k\e,\e)=\exp\left(\e^{-1} \Phi_1(k\e,\e)\right),
$$ 
\item[(b)] $\Phi_1$ has an asymptotic expansion (uniform with respect
to $x$):
$$
\Phi_1(x,\e) \sim_{\e\to 0} \exp\left( \e^{-1}
\sum_{s=0}^\infty \phi_{1,s}(x) \e^s \right)
$$
where $\phi_{1,s}$ are as in Lemma \ref{lem.borel}. As a result,
we have an asymptotic expansion (uniform with respect to $k$):
$$
\psi_{1}(k\e,\e) \sim_{\e \to 0} \exp\left(\e^{-1} \sum_{s=0}^\infty 
\phi_{1,s}(k\e)  \e^s\right).
$$
\end{itemize}
\end{proposition}

\begin{proof}
Consider the change of variables $\th$ as in \eqref{eq.change}.

Due to our choice of initial conditions it follows that for every
fixed $k=0,\dots,d-1$, the function $\e\to\th(k\e,\e)$ is smooth. Using this
and the smoothness of the coefficients of \eqref{eq.ediff}, it follows that
for every fixed $k$, the function $\e\to\th(k\e,\e)$ is smooth.

So far, the function $\th$ is defined on $\D^{\mathrm{seg}}_{\e',I}$:
$$
\psdraw{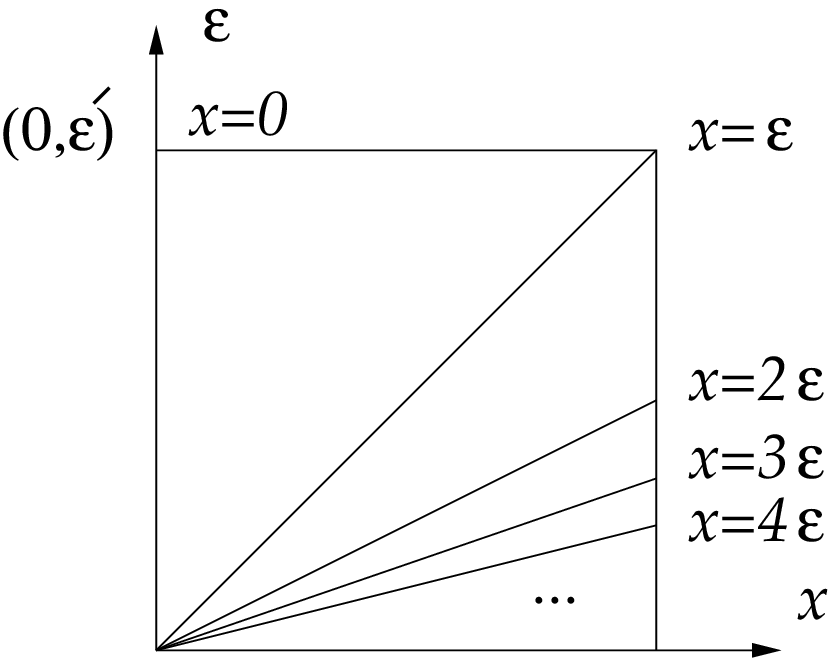}{2in}
$$
which is a union of line segments in a rectangle $I \times [0,\e']$, and 
satisfies \eqref{eq.estimate}.

The complement of these line segments in $[0,\e']$
consists of an infinite union of open triangles, together with the horizontal
segment $I \times 0$. 
We can smoothly interpolate $\th$ inside these open triangles so that
it is defined on $I \times (0,\e']$ and 
\begin{equation}
\lbl{eq.intfinal}
|\th(x,\e)-1| < C_s \e^s
\end{equation}
for all $(x,\e) \in I \times (0,\e']$ and all $s \in \BN$.

Let us extend $\th$ to $I \times [0,\e']$ by defining $\th(x,0)=1$
for all $x \in I$.

We claim that $\th$ is smooth on $I \times [0,\e']$. We need only to check
this at the points $(x,0)$ for $x \in I$. This follows easily from 
\eqref{eq.intfinal}. For example, to check continuity at $(x,0)$, consider
a sequence $(x_n,\e_n)$ such that $\lim_{n\to\infty} (x_n,\e_n)=(x,0)$.
Then, \eqref{eq.intfinal} for $s=1$ implies that
$|\th(x_n,\e_n)-1| < C_s \e_n \sim_{n\to\infty} 0$, thus $\th$ is
continuous at $(x,0)$. Using \eqref{eq.intfinal} for $s+1$ it follows that
$\partial^{s}/\partial \e^s \th(x,\e)|_{\e=0}=0$ for all $s>0$, and we find
that $\th$ has an $\e$-asymptotic expansion (uniform with respect to $x$):
$$
\th(x,\e) \sim_{\e\to 0} 1
$$
Restricting further $\e'$ if needed, we may assume that $|\th(x,\e)| >0$
for all $(x,\e) \in I \times [0,\e']$; in other words $\log\th(x,\e)$
makes sense for all $(x,\e) \in I \times [0,\e']$. 

Now, we can finish the proof of Proposition \ref{prop.step1}.

Let us define
$$
\Phi_1(x,\e)=\hat\Phi_1(x,\e)+\e \log \th(x,\e)
$$
Then, \eqref{eq.change} implies (a). Since $\th(x,\e)$ is asymptotic to
$1$ (uniformly on $x$), it follows that $\Phi_1(x,\e)$ is asymptotic
to $\hat\Phi_1(x,\e)$. Using the asymptotic of $\hat\Phi_1$ given by
Lemma \ref{lem.borel}, (b) follows.
\end{proof}

\subsection{A reduction to an $\e$-difference equation of smaller degree}
\lbl{sub.reduction}

We will now prove Theorem \ref{thm.2} by induction on the degree $d$
of the $\e$-difference equation. For $d=1$, it follows from Proposition
\ref{prop.step1}. The inductive step is the next Proposition.

\begin{proposition}
\lbl{prop.step2}
Assume that Theorem \ref{thm.2} holds for regular $\e$-difference
equations of degree less than $d$. Then it holds for regular
$\e$-difference equations of degree $d$. 
\end{proposition}

\noindent
{\em Proof.}
Consider a $\e$-difference equation \eqref{eq.ediff} of degree $d$.
We will use the solution $\psi_1$ of it constructed in Proposition
\ref{prop.step1} to reduce it to an equivalent equation of degree $d-1$,
and an inhomogeneous $\e$-difference equation of degree $1$.

Consider the dependent change of variables 
\begin{equation}
\lbl{eq.dependent}
\th=\frac{\phi}{\psi_1}
\end{equation}
This is well-defined since $\psi_1$ is nowhere zero. Then, $\phi$ satisfies
\eqref{eq.ediff} iff $\th$ satisfies 
\begin{equation}
\lbl{eq.reduce}
\sum_{j=0}^d b_j(k\e,\e) \th((k+j)\e,\e)=0
\end{equation}
where
$$
b_j(x,\e)=a_j(x,\e) \frac{\psi_1(x+j\e,\e)}{\psi_1(x,\e)}.
$$ 
The characteristic polynomials of \eqref{eq.ediff} and \eqref{eq.reduce}
are related by
$$
P_{\eqref{eq.reduce}}(\l)=\l_1(x)^d P_{\eqref{eq.ediff}}(\l/\l_1(x)).
$$
As in the proof of Proposition \ref{prop.step1}, it is easy to see that
\eqref{eq.reduce} is a regular $\e$-difference equation.
Moreover, it is easy to see that Theorem \ref{thm.2} holds for \eqref{eq.ediff}
iff it holds for \eqref{eq.reduce}. Indeed, check that the change of
variables given by \eqref{eq.dependent} preserves the asymptotics of the
solutions of \eqref{eq.ediff} and \eqref{eq.reduce}.

Thus, it suffices to work with \eqref{eq.reduce}. In that case, 
$\th=1$ is a solution of \eqref{eq.reduce},
since $\psi$ is a solution of \eqref{eq.ediff}. It follows that
\begin{equation}
\lbl{eq.zz}
\sum_{j=0}^d b_j(x,\e)=0.
\end{equation}
(Compare this with \eqref{eq.approximate}).
Let us define 
\begin{equation}
\lbl{eq.zeta}
\zeta(k\e,\e)=\th((k+1)\e,\e)-\th(k\e,\e).
\end{equation}
Then, we get that $\zeta$ is a solution of the $\e$-difference equation
\begin{equation}
\lbl{eq.reduce1}
\sum_{j=0}^{d-1} c_j(k\e,\e) \zeta((k+j)\e,\e)=0
\end{equation}
where
$$
c_s(x,\e)=\sum_{j=s+1}^d b_j(x,\e).
$$
The characteristic equations of \eqref{eq.reduce} and \eqref{eq.reduce1}
are related by
$$
\sum_{j=0}^d b_j(x,0) \l^j= (\l-1) \sum_{j=0}^{d-1} c_j(x,0) \l^j
$$
Since $c_0(x,\e)=\sum_{j=1}^d b_j(x,\e)=-b_0(x,\e)$ (by \eqref{eq.zz})
and $c_d(x,\e)=b_d(x,\e)$, the same arguments of Proposition \ref{prop.step1} 
imply that \eqref{eq.reduce1} is regular, 
assuming that \eqref{eq.reduce} is regular.

By the induction hypothesis, it follows that \eqref{eq.reduce1} satisfies
Theorem \ref{thm.2}. 

For the remainder of this section, fix 
a solution $\zeta$ of \eqref{eq.reduce1} which satisfies the
properties of Theorem \ref{thm.2}. In other words, $\zeta$ satisfies 
\eqref{eq.reduce} and $\zeta(k\e,\e)=\exp(\e^{-1} Z(k\e,\e))$ where $Z$ is
a smooth function with uniform (with respect to $x$) asymptotics:
$$
Z(x,\e)\sim_{\e\to 0} \sum_{s=0}^\infty Z_s(x) \e^s.
$$

\begin{lemma}
\lbl{lem.stepa}
There exists a formal solution 
$$
\tilde{\th}(x,\e)=\exp(\e^{-1} \Th(x,\e))
$$ 
of \eqref{eq.reduce1} such that
$$
\Th(x,\e)=\sum_{s=0}^\infty \Th_s(x) \e^s.
$$
\end{lemma}

\begin{proof}
We need to solve the formal power series Equation
$$
\exp\left(\e^{-1} \sum_{s=0}^\infty \Th_s(x+\e)\e^s\right)
-
\exp\left(\e^{-1} \sum_{s=0}^\infty \Th_s(x)\e^s\right)=
\exp\left(\e^{-1} \sum_{s=0}^\infty Z_s(x)\e^s\right)
$$
for $\Th$ in terms of $Z$.
Using the Taylor expansion $\Th_0(x+\e)=\Th_0(x)+\Th'_0(x)\e + O(\e^2)$
it is easy to see that the above equation
equals to
$$
\exp\left(\e^{-1} \Th_0(x) + O(1)\right) 
-
\exp\left(\e^{-1} \Th_0(x) + O(1) \right)=
\exp\left(\e^{-1} Z_0(x) + O(1) \right)
$$
from which follows that $\Th_0=Z_0$. Dividing the equation by
$\exp(\e^{-1}\Th_0)$ we get an equation in formal power series with nonnegative
powers of $\e$.
Moreover, the coefficient of $\e^s$ in that power series (for $s \geq 0$)
equals to
$$
(\exp(\Th_0'(x))-1) \exp(\Th_{s+1}(x)) H_s(x)
$$
where $H_s(x)$ is a function of $\Th_j$ and $Z_j$ for $j=1,\dots,s$.
Since $\exp(\Th_0'(x))=\exp(Z_0'(x))$ is an eigenvalue of \eqref{eq.reduce1},
it is never equal to $1$.

This and induction prove  the lemma. 
\end{proof}

\begin{lemma}
\lbl{lem.stepb}
\rm{(a)} There exists a solution to the equation
\begin{equation}
\lbl{eq.zeta2}
\zeta(k\e,\e)=\theta((k+1)\e,\e)-\theta(k\e,\e)
\end{equation}
for $\theta$ in terms of 
$\zeta$ with appropriate initial condition.
\newline
\rm{(b)} For all $(k\e,\e) \in \D'_{\e',I}$ we have:
$$
\zeta(k\e,\e) \sim_{\e\to 0} \exp(Z(k\e,\e))
$$
where 
$$
\exp(Z(x,\e))=\exp(\Th(x+\e,\e))-\exp(\Th(x,\e))
$$
\end{lemma}

\begin{proof}
Fix $\e > 0$ and let $I=[a,b]$. Consider a natural number $k$ such that
$k\e \in I$ and $(k+1)\e \in I$. This is equivalent to 
$k_1 \leq k \leq k_2$ where $k_1$ and $k_2$ are natural numbers 
that depend on $\e$ and $I$,
although we do not indicate this in our notation.

Then equation \eqref{eq.zeta2} implies that
\begin{eqnarray*}
\zeta(k_1\e,\e) &=& \th((k_1+1)\e,\e)-\th(k_1\e,\e) \\
\zeta((k_1+1)\e,\e) &=& \th((k_1+2)\e,\e)-\th((k_1+1)\e,\e) \\
\dots &=& \dots \\
\zeta((k_2-1)\e,\e) &=& \th(k_2\e,\e)-\th((k_2-1)\e,\e) 
\end{eqnarray*}
Summing up, we obtain that
$$
\th(k_2\e,\e)=\th(k_1\e,\e)+\sum_{j=k_1}^{k_2-1} \zeta(j\e,\e).
$$
Choose initial conditions so that $\th(k_1\e,\e)=\exp(Z(k_1\e,\e))$.
This completes part (a).

Part (b) follows by a telescoping calculation.
\end{proof}

To finish the proof of Proposition \ref{prop.step2}, it suffices to show
that the solution $\zeta$ of Lemma \ref{lem.stepb} is asymptotic to the
formal solution $\tilde{\th}$ of Lemma \ref{lem.stepa}.

Since 
$$
\zeta(k\e,\e) \sim_{\e\to 0} \exp(Z(k\e,\e))
$$ 
and 
$$
\exp(Z(k\e,\e)) = \exp(\Th((k+1)\e,\e))-\exp(\Th(k\e,\e)) 
$$
it follows by the definition of $\th$ given in Lemma \ref{lem.stepb} 
and by a telescopic sum, that:

\begin{eqnarray*}
\th(k\e,\e) &=& \th(k_1\e,\e) +\sum_{j=k_1}^{k-1} \zeta(j\e,\e) \\
& \sim_{\e\to 0} &  \th(k_1\e,\e) +\sum_{j=1}^{k-1} 
\left(\exp(\Th((j+1)\e,\e))-\exp(\Th(j\e,\e)) \right) \\
&=& \th(k_1\e,\e) + \exp(\Th(k\e,\e)) - \exp(\Th(k\e,\e)) \\
&=& \exp(\Th(k\e,\e))
\end{eqnarray*}

This concludes the proof of Proposition \ref{prop.step2}.
\qed

\subsection{The solutions form a locally fundamental set}
\lbl{sub.form}

Let us summarize what we have obtained so far.

Consider a partition $I=\cup_{p \in P} I_p$ 
of $I=[x_0,x_P]$ given by
$I_p=[x_p, x_{p+1}]$ for $p=1,\dots,P-1$, and consider a permutation
$\s_p$ of $\{1,\dots,d\}$ such that
$$
|\l_{\s_p(1)}(x)| \geq |\l_{\s_p(2)}(x)| \geq \dots \geq
|\l_{\s_p(d)}(x)| \qquad \text{for all} \qquad x \in I_p. 
$$
We have constructed 
solutions smooth functions $\Phi_m$ for $m=1,\dots,d$ with
asymptotic expansion given by \eqref{eq.Phim} and \eqref{eq.intc}.

Let us define
\begin{equation}
\lbl{eq.fornow}
\psi_m(x,\e)=\exp\left(\e^{-1} \Phi_m(x,\e)\right),
\end{equation}
where $\Phi_m$ are smooth functions with asymptotic expansions
as in Equations \eqref{eq.Phim} and \eqref{eq.intc}.

Moreover, we have shown that for every interval $I_p$ (as in the discussion
prior to Theorem \ref{thm.2}), 
$$
\{\psi_1(k\e,\e),\dots, \psi_d(k\e,\e)\}
$$
is a set of solutions of \eqref{eq.ediff} when $k\e \in I_p$.

Fix a solution $\psi(k\e,\e)$ of \eqref{eq.ediff} and 
an interval $I_p$. The following lemma certainly implies that $\{\psi_1,\dots,
\psi_d\}$ is a locally fundamental set of solutions of \eqref{eq.ediff}.
This concludes the proof of Theorem \ref{thm.2}.
\qed

In addition, the next lemma
 motivates the definition of a regular solution, given in the
following section.

\begin{lemma}
\lbl{lem.certainly}
\rm{(a)} Fix $\psi$ and $I_p$ as above. Then, there exist smooth functions
$c^p_m$ such that
\begin{equation}
\lbl{eq.psinow}
\psi(k\e,\e)=c^p_1(\e) \psi_{\s_p(1)}(k\e,\e) + \dots + 
c^p_d(\e) \psi_{\s_p(d)}(k\e,\e)
\end{equation}
for all $k\e \in I_p$. \newline
\rm{(b)}
Moreover, for every $p$ and $m=1,\dots,d$ we have
\begin{equation}
\lbl{eq.cpm}
c^p_m(\e)=\frac{\psi_{\s_{p-1}(m)}(\xpe\e-\e,\e)}{
\psi_{\s_{p}(m)}(\xpe\e,\e)}
\ga^p_m(\e)
\end{equation}
for some smooth functions $\ga^p_m$, with the understanding that  
$\psi_{\s_{-1}(m)}=1$. Here $[x]$ is the {\em largest integer smaller than} 
$x$, and
$$
\xpe=\left[\frac{x_p}{\e}\right] +1.
$$
\end{lemma}

\begin{proof}
Without loss of generality, let us assume that $\s_p(j)=j$ for $j=1,\dots,d$.
Equation \eqref{eq.psinow} is a linear equation in $c^p_m$, with solutions
$$
c^p_m(\e)=\frac{\det
W_m(\xpe,\e)}{\det W(\xpe,\e)}
$$
where $I_p=[x_p,x_{p+1}]$, 
$$
W(x,\e)=
\begin{pmatrix}
\psi_1(x,\e) & \dots & \psi_m(x,\e) & \dots & \psi_d(x,\e) \\
\psi_1(x+\e,\e) & \dots & \psi_m(x+\e,\e) & \dots & \psi_d(x+\e,\e) \\
\dots & \dots & \dots & \dots & \dots \\
\psi_1(x+(d-1)\e,\e) & \dots & \psi_m(x+(d-1)\e,\e) & \dots &
\psi_d(x+(d-1)\e,\e) \\
\end{pmatrix}
$$
and 
$$
W_m(x,\e)=
\begin{pmatrix}
\psi_1(x,\e) & \dots & \psi(x,\e) & \dots & \psi_d(x,\e) \\
\psi_1(x+\e,\e) & \dots & \psi(x+\e,\e) & \dots & \psi_d(x+\e,\e) \\
\dots & \dots & \dots & \dots & \dots \\
\psi_1(x+(d-1)\e,\e) & \dots & \psi(x+(d-1)\e,\e) & \dots & 
\psi_d(x+(d-1)\e,\e) \\
\end{pmatrix}
$$
where the $\psi$'s are in the $m$th column of $W_m$.

We will show that for small enough $\e$, $W(x,\e)$ is nonsingular.

Using Equations \eqref{eq.fornow}, \eqref{eq.Phim} and \eqref{eq.intc},
it follows that
 
\begin{eqnarray*}
\frac{1}{\psi_1(x,\e) \dots \psi_d(x,\e)} \det(W(x,\e)) &=&
\det( \exp(j \phi'_m(x)) ) + O(\e) \\ 
 &=& \det (\l_m(x)^j) + O(\e) \\
&=& \pm \prod_{i \neq j}(\l_i(x)-\l_j(x)) + O(\e)
\end{eqnarray*}
uniformly in $x$, where $\l_m$ are the eigenvalues of \eqref{eq.ediff}.
Since \eqref{eq.ediff} is regular, its eigenvalues never collide.

Similarly, using Equation \eqref{eq.fornow}, we have:
$$
\frac{1}{\psi_1(\xpe\e,\e) \dots \widehat{\psi_m}(\xpe\e,\e) 
\dots \psi_d(\xpe\e,\e)} 
\det(W_m(\xpe\e,\e)) 
=
\det B_{m,0}(\xpe\e,\e) + O(\e)
$$
where
$$
B_{m,j}(x,\e)=
\begin{pmatrix}
\l_1(x)^{-j} & \dots & \psi(x,\e) & \dots & \l_d(x)^{-j} \\
\l_1(x)^{1-j} & \dots &  \psi(x+\e,\e) & \dots & \l_d(x)^{1-j} \\
\dots & \dots & \dots & \dots & \dots \\
\l_1(x)^{d-1-j} & \dots & \psi(x+(d-1)\e,\e) & \dots & 
\l_d(x)^{d-1-j} \\
\end{pmatrix}
$$
where the $\psi$'s are in the $m$th column of $B_{m,j}$.
Thus,
\begin{equation}
\lbl{eq.tmpg}
c^p_m(\e)= \frac{1}{\psi_m(\xpe\e,\e)} 
\frac{\det(B_{m,0}(\xpe\e,\e))}{\prod_{j \neq m} \l_j(\xpe\e)-\l_m(\xpe\e)}
+ O(\e).
\end{equation}

The idea now is to move the recursion relation
backwards $d$ times. Using the solution $\psi(k\e,\e)$ for $k\e \in I_{p-1}$
will allow us to compute the smooth functions $c^p_m$.

In detail, consider the matrix $B_{m,0}(\xpe\e,\e)$ and move the recursion 
relation backwards once. 
Using Equation \eqref{eq.matrixform} and the fact that the 
the $j$th column of $B_{m,j}(\xpe,\e)$ for $j \neq m$ is an eigenvector of
$A(\xpe\e,\e)$ (up to $O(\e)$) with eigenvalue $\l_j(\xpe\e)$, it follows that
$$
B_{m,0}(\xpe\e,\e)=A(\xpe\e,\e) B_{m,1}(\xpe\e-\e,\e) + O(\e).
$$
Iterating $d-1$ more times, it follows that
$$ 
B_{m,0}(\xpe\e,\e)=A(\xpe\e,\e) \dots A(\xpe\e-(d-1)\e,\e)
B_{m,d}(\xpe\e-d\e,\e)  + O(\e).
$$
Since $\xpe\e-d\e \in I_{p-1}$, it follows (as in the computation of 
$W_m(x,\e)$ above) that:
\begin{eqnarray*}
\det(B_m(\xpe\e-d\e,\e)) &=& \psi_{\s_{p-1}(m)}(\xpe\e-\e,\e) \\
& \det &  \begin{pmatrix}
\l_1(\xpe\e-\e)^{1-d} & \dots & \frac{\psi_{\s_{p-1}(m)}(\xpe\e-d\e,\e)}{
\psi_{\s_{p-1}(m)}(\xpe\e-\e,\e)} & \dots & \l_d(\xpe\e-\e)^{1-d} \\
\l_1(\xpe\e-\e)^{2-d} & \dots & \frac{\psi_{\s_{p-1}(m)}(\xpe\e+(1-d)\e,
\e)}{
\psi_{\s_{p-1}(m)}(\xpe\e-\e,\e)} & \dots & \l_d(\xpe\e-\e)^{2-d} \\
\dots & \dots & \dots & \dots & \dots \\
\l_1(\xpe\e-\e)^{d-d} & \dots & \frac{\psi_{\s_{p-1}(m)}(\xpe\e-\e,
\e)}{
\psi_{\s_{p-1}(m)}(\xpe\e-\e,\e)} & \dots & 
\l_d(\xpe\e-\e)^{d-d} \\
\end{pmatrix}
\\ &=& 
\psi_{\s_{p-1}(m)}(\xpe\e-\e,\e) \\
& & \prod_{i \neq j, j
\neq m}
(\l_i^{-1}(\xpe\e-\e)-\l_j^{-1}(\xpe\e-\e)) \\ & & 
\prod_{i \neq m}
(\l_i^{-1}(\xpe\e-\e)-\l_{\s_{p-1}(m)}^{-1}(\xpe\e-\e)) +O(\e).
\end{eqnarray*}
This, together with Equation \eqref{eq.tmpg} proves \eqref{eq.cpm}.
\end{proof}

\section{Regular solutions and their asymptotics}
\lbl{sec.asymptotics}

In this section we discuss regular solutions of $q$ and $\e$-difference
equations and their asymptotics.

\subsection{Regular solutions to $\e$-difference equations}
\lbl{sub.regulare}

According to Lemma \ref{lem.certainly}, a solution $\psi$ to \eqref{eq.ediff} 
determines 
a collection $S=\{S_p|\, p \in \P\}$ of subsets of $\{1,\dots, d\}$, where
$$
S_p=\{m \in \{1,\dots,d\} | \ga^p_m \neq 0\}.
$$

\begin{definition}
\lbl{def.regularsol}
Fix a collection $S=\{S_p|\, p \in \P\}$ of subsets of $\{1,\dots, d\}$.
We say that a solution $\psi$ of a regular equation 
\eqref{eq.ediff} is $S$-{\em regular}
iff for every $p \in \P$ we have:
\begin{itemize}
\item
$\ga^p_m=0$ if $m \not\in S_p$.
\item
$\ga^p_m$ are {\em not flat} at $0$ for all $m \in S_p$. 
That is, some derivative of $\ga^p_m$ at $\e=0$ does not vanish.
\item
For every $p \in \P$ there exists an element $\eta(p) \in S_p$ such that
$$
|\l_{\eta(p)}(x)| > |\l_j(x)|  \quad \text{for all} \quad j \in 
S_p-\{\eta(p)\}, \quad 
x \in \text{Interior}(I_p).
$$ 
In other words, in the interior of the interval $I_p$, and among
the eigenvalues $\l_j(x)$ for $j \in S_p$, there is a unique
eigenvalue of strictly maximum magnitude.
\end{itemize}
We will say that a solution to \eqref{eq.qdiff} is {\em regular} if it
is $S$-regular for some $S$.
\end{definition}

\subsection{Asymptotics of regular solutions of $\e$-difference equations}
\lbl{sub.assregulare}

\begin{proof}(of Theorem \ref{thm.ass})
Let $\psi$ be an $S$-regular solution to \eqref{eq.ediff}. 
Let us assume that $S=\{1,\dots,d\}$, and  
$|\l_1(x)| >
|\l_m(x)|$ for $m=2, \dots, d$ and all $x$ in the interior $\text{Interior}(I)$
of the closed interval $I=[0,b]$. Fix an $x \in I$.

Then, we have:
$$
\psi(x,\e)=c_1(\e) \psi_1(x,\e) + \dots + c_d(\e) \psi_d(x,\e).
$$  
for $x=k\e$, where $c_1(\e)=c_1 \e^{n_1} + O(\e^{n_1+1})$, and $c_1 \neq 0$.

Then, we have:
\begin{equation}
\lbl{eq.now}
\psi(x,\e) = c_1(\e) \psi_1(x,\e) \left(
1+ \sum_{m=2}^d \frac{c_j(\e)}{c_1(\e)} \frac{\psi_m(x,\e)}{\psi_1(x,\e)}
\right).
\end{equation}
Recall from Theorem \ref{thm.2} that
\begin{eqnarray*}
\psi_m(x,\e) &=& \exp(\e^{-1} \Phi_m(x,\e)) \\
\Phi_m(x,\e) &=& \Phi_{m,0}(x) + O(\e) \\
\Phi_{m,0}'(x) &=& \log\l_m(x).
\end{eqnarray*}
Thus, 
$$
\text{Re}(\Phi_{m,0})'(x)=\text{Re}(\log\l_m(x))=\log|\l_m(x)|.
$$ 
Combined with $|\l_m(x)| < |\l_1(x)|$
for $m \geq 2$, and $\Phi_{m,0}(0)=0=\Phi_{1,0}(0)$,
it follows that 
$$
\text{Re}(\Phi_{m,0})(x) < \text{Re}(\Phi_{1,0})(x)
$$
for all $x \in I$ and
$$
\text{Re}(\Phi_{1,0})(x)=\int_0^x \log| \l_1(t)| dt.
$$
Therefore, 
$$
\lim_{\e\to 0^+}
\left(
1+ \sum_{m=2}^d \frac{c_j(\e)}{c_1(\e)} \frac{\psi_m(x,\e)}{\psi_1(x,\e)}
\right)=1.
$$
and
$$
\lim_{\e\to 0^+} \e\log|c_1(\e)|=0.
$$
Thus, Equation \eqref{eq.now} implies that
\begin{eqnarray*}
\lim_{\e\to 0^+} \e \log |\psi(x,\e)| &=&
\lim_{\e\to 0^+} \e \log |\psi_1(x,\e)| \\ &=&
\lim_{\e\to 0^+} \e \log |\exp(\e^{-1} \Phi_m(x,\e))| \\ & = &
\lim_{\e\to 0^+} \e \log |\exp(\e^{-1} \Phi_m(x,0) + O(1))| \\ &=&
\text{Re}(\Phi_{1,0}(x)) \\ &=& \int_0^x 
\log| \l_1(t)| dt
\end{eqnarray*}
The result follows. 

In the general case, we partition the interval $I_p$ as in the discussion
prior to Theorem \ref{thm.2} and repeat the above proof using 
Equation \eqref{eq.gpm}. The result follows. 
\end{proof}

\subsection{Asymptotics of regular solutions of $q$-difference equations}
\lbl{sub.assregularq}

First, we need to define what is a regular solution to a $q$-difference
equation.

Consider a solution $\psi$ of a $q$-difference equation and a partition
of $S^1$ as in Section \ref{sub.qdiff}. Then, at each interval $I_p$,
we can write the solution as a linear combination of fundamental solutions,
as in Equation \eqref{eq.cpms}. Let $S_p$ be the indexing set of the
fundamental solutions that we use in each interval $I_p$; see \eqref{eq.Sp}.

Suppose
that the partition of $S^1$ is given by
$I_p=[e^{2\pi i x_p}, e^{2\pi i x_{p+1}}]$ for $p=0,\dots,P-1$, with
$q_0=1,\ q_P=e^{2\pi i}$. 

Then, with $q=e^{2\pi i \e}$ it turns out that for every 
$p$ and $m=1,\dots,d$ we have
\begin{equation}
\lbl{eq.gpm}
c^p_m(q)=\frac{\psi_{\s_{p-1}(m)}(\xpe\e-\e,\e)}{\psi_{\s_{p}(m)}(\xpe\e,\e)}
\ga^p_m(q)
\end{equation}
for some smooth functions $\ga^p_m$, with the understanding that  
$\psi_{\s_{-1}(m)}=1$.

\begin{definition}
\lbl{def.regularsolq}
Fix a collection $S=\{S_p|\, p \in \P\}$ of subsets of $\{1,\dots, d\}$.
We say that a solution $\psi$ of a regular equation 
\eqref{eq.qdiff} is $S$-{\em regular}
iff for every $p \in \P$ we have:
\begin{itemize}
\item
$\ga^p_m=0$ if $m \not\in S_p$.
\item
$\ga^p_m$ are {\em not flat} at $0$ for all $m \in S_p$. 
That is, some derivative of $\ga^p_m$ at $\e=0$ does not vanish.
\item
For every $p \in \P$ there exists an element $\eta(p) \in S_p$ such that
$$
|\l_{\ga(p)}(v)| > |\l_j(v)|  \quad \text{for all} \quad j \in 
S_p-\{\ga(p)\}, \quad 
v \in \text{Interior}(I_p).
$$ 
In other words, in the interior of the interval $I_p$, and among
the eigenvalues $\l_j(v)$ for $j \in S_p$, there is a unique
eigenvalue of strictly maximum magnitude.
\end{itemize}
We will say that a solution to \eqref{eq.qdiff} is {\em regular} if it
is $S$-regular for some $S$.
\end{definition}

\begin{proof}(of Theorem \ref{thm.assq})
It follows from Equation \eqref{eq.phiff} of Lemma \ref{lem.translate} 
and Theorem \ref{thm.ass}.
\end{proof}

\section{Applications to Quantum Topology}
\lbl{sub.goodbad}

\subsection{The $A$-polynomial of a knot and its noncommutative version}
\lbl{sub.A}

In this section we discuss general features of $q$-difference equations
for the colored Jones function of a knot.

The coefficients of the $q$-difference equations are rational functions
of $q$ and $q^n=Q$. In order to simplify the typesetting, we will give the
$q$-difference equation 
$$
\sum_{j=0}^d b_j(q^n,q) f(n+j)=0
$$
in operator form 
\begin{equation}
\lbl{eq.operator}
Pf=0
\end{equation}
where
$$
P=\sum_{j=0}^d b_j(Q,q)E^j
$$
and 
where the operators $E$, $Q$ and $q$,
act on a discrete function $f:\BN \to \BZ[q^{\pm}]$ by
$$
(q f)(n)=q f(n) \qquad (Qf)(n)=q^n f(n) \qquad (Ef)(n)=f(n+1).
$$
Note that $q$ commutes with $Q$ and $E$, and that $EQ=qQE$.

It follows by definition that the characteristic polynomial $\ch P(v,\l)$ of 
\eqref{eq.operator} is obtained from $P$ by setting $q=1$, 
replacing $(E,Q)$ by $(\l,v)$. In other words, we have:
$$
\ch P(v,\l)=\sum_{j=0}^d b_j(v,1)\l^j.
$$

In \cite{Ga1}, the first author showed that the colored Jones function $J_K$
of a knot $K$ satisfies an essentially unique smallest degree $q$-difference
equation $P_K J_K=0$ where the coefficients $a_j(u,v)$ of $P_K$ are rational
functions of $u$ and $v$ with rational coefficients.

In \cite{Ga1}, the operator $P_K$ was called the 
{\em non-commutative $A$-polynomial} of $K$. 

In \cite{Ga1}, it was conjectured that:
 
\begin{conjecture}(AJ Conjecture)
Up to a multiplication by a polynomial in $v$, we have
$$
\ch P_K(\l,v)=A_K(L,M)|_{(L,M^2)=(\l,v)}
$$
where $A_K$ is the $A$-polynomial of $K$, defined by \cite{CCGLS}.
\end{conjecture}

The $A$-polynomial of $K$ parametrizes the
moduli space of characters of $\SL_2(\BC)$ representations of $\pi_1(S^3-K)$,
restricted to the boundary torus $\pt M$. The $A$-polynomial of a knot
is an important ingredient to the Geometrization of the knot complement 
and its Dehn fillings.

The $A$-polynomial $A_K$ of a knot $K$ in $S^3$ satisfies symmetries, which
we will list here, and refer to \cite{CCGLS} and \cite{CL} for proofs.

\begin{itemize}
\item[{\bf (S1)}]
It has integer coefficients and even powers of $M$, that is
 $A_K(L,M) \in \BZ[L,M^2]$.
\item[{\bf (S2)}]
It is reciprocal, that is, 
$A_K(L^{-1},M^{-1})=\pm L^k M^l A(L,M)$.
\item[{\bf (S3)}]
It is tempered, that is the edge polynomials of its Newton polygon
are cyclotomic.
\item[{\bf (S4)}]
It specializes to 
$$
A_K(L,1)=\pm (L-1)^{n_+} (L+1)^{n_-}
$$
for some integers $n_{\pm}$.
\item[{\bf (S5)}]
$L-1$ is always a factor of $A_K$, that corresponds to $U(1)$ representations.
\end{itemize}

If the colored Jones function of a knot was an $S$-regular solution to a 
regular
$q$-difference equation, and if the AJ Conjecture were true, 
then it follows that for every $\a \in [0,1]$ we have:
\begin{equation}
\lbl{eq.HVCA}
\lim_{n \to \infty} \frac{\log|J_K(n)(e^{\frac{2 \pi i \a}{n}})|}{n}
=
\s^J_{S,K}(\a)=\s^A_{S,K}(\a)
\end{equation}
where $\s^A_{S,K}$ is the {\em $A$-entropy} 
of a knot, defined as follows.

\begin{definition}
\lbl{def.Aentropy}
For a knot $K$ in $S^3$, let $L_j(t)$ for $j=1,\dots, d$ denote the roots
of the equation
$$
A_K(L_j(t),e^{i t/2})=0
$$
for $t \in [0,2 \pi]$, where $d$ is the $L$-degree of $A_K$. 
Fix a partition $\cup_{p\in \P}I_p$ of $[0,2 \pi]$
by closed intervals with nonoverlapping interiors and a permutation
$\s_p$ of the set $\{1,\dots,d\}$ such that
$$
|L_{\s_p(1)}(t)| \geq |L_{\s_p(2)}(t)| \geq \dots \geq
|L_{\s_p(d)}(t)| \qquad \text{for all} \qquad t \in I_p. 
$$
For every collection $S=\{S_p|\, p \in \P\}$ of subsets of $\{1,\dots, d\}$,
we can define the $S$-{\em entropy} 
$$
\s^A_{S,K}:[0,1]\to\BR$$
by
$$
\s^A_{S,K}(\a)=\frac{1}{2 \pi}
\int_0^{2 \pi} \log \chi_S(\a t) dt,
$$
where $\chi_S: [0,1] \to \BR$ is defined by 
$$
\chi_S(t)=\max_{j \in S_p} |L_{\s_p(j)}(t)| \qquad \text{if} \qquad
t \in I_p.
$$
\end{definition}

It is natural to ask how the $A$-entropy of a hyperbolic knot (evaluated
at $\a=1$) compares to
the Hyperbolic Volume. The answer to this question is essentially contained
in work of D. Boyd, \cite{Bo}, which we quote without proof here.
We urge the reader to look in \cite{Bo} for beautiful and suggestive
calculations.
 
Boyd introduced and studied 
another invariant of a knot, the {\em Mahler measure}
$$
m_K=\int_{S^1 \times S^1} \log|A_K(x,y)| dx dy
$$
Using Jensen's formula, and the symmetry (S3), it follows that
$$
2 \pi m_K=\sum_{j=0}^d  \int_0^{2 \pi} \log^+|L_j(t)| dt
$$
where $\log^+(a)=\log \max\{1,a\}$.

Using (S2), it follows that
$1/L_j(t)$ is an eigenvalue for every eigenvalue $L_j(t)$. 

More generally, among the roots $L_j(t)$ there is a distinguished one, $L_1$,
corresponding to the discrete faithful representation when $t=2\pi$. 
Let $L_d(t)=1$
denote the eigenvalue corresponding to the $U(1)$ representations.
Boyd informs us that
for $2$-bridge knots $K$ (in particular, for the $3_1$ and $4_1$ knots), it
is true that
$$
\s^A_{\{1,d\},K}(1)=\text{vol}_K.
$$

It follows by (S4) that the eigenvalues collide at $t=0$. 
Moreover, Boyd informs us that for $2$-bridge knots there exists a $t_0 \in
(0, \pi)$ such that $|L_j(t)|=1$ for all $j$ and all $t: t_0 < t <  \pi$.

Thus, if we want to apply Theorem \ref{thm.1} to the GHVC, we need to
deal with irregular $q$-difference equations. We will discuss this
topic in detail in a later publication. 

Meanwhile, let us discuss some examples, taken from \cite{Ga1}.

\subsection{Examples: The $3_1$ and $4_1$ knots}
\lbl{sub.examples}

In this section we discuss in detail $q$-difference equation of the 
the colored Jones function of the
two simplest knots, namely the trefoil $3_1$ and the figure eight $4_1$.
The former is not hyperbolic, and the latter is.

In \cite{Ga1}, the first author 
computed that the colored Jones function $J_{3_1}$
(resp. $J_{4_1}$) satisfies the second (resp. third) order $q$-difference
equation
$$
P_{3_1} J_{3_1}=0 \qquad \text{resp.} \qquad 
P_{4_1} J_{4_1}=0
$$
where the noncommutative $A$-polynomials $P_{3_1}$ and $P_{4_1}$ are given
by:

\begin{eqnarray*}
P_{3_1}&=&
\frac{q^3 Q^2 (q^2 - q^2 Q)}{q^3 - q^4 Q^2} \\ & & + \frac{(q - 
q^2 Q) (q + q^2 Q) (q^4 - q^5 Q + q^6 Q^2 - q^7 Q^2 - 
q^7 Q^3 + q^8 Q^4)}{q^2 Q (q - q^4 Q^2) (q^3 - q^4 Q^2)}E  \\
& & + \frac{-1 + q^2 Q}{Q (q - q^4 Q^2)}E^2
\\
P_{4_1}&=&
\frac{q^5 Q (-q^3 + q^3 Q)}{(q^2 + q^3 Q) 
(-q^5 + q^6 Q^2)} \\ 
& & - \frac{(q^2 - q^3 Q) (q^8 - 
              2 q^9 Q + q^{10} Q - q^9 Q^2 + q^{10} Q^2 - 
              q^{11} Q^2 + q^{10} Q^3 - 2 q^{11} Q^3 + 
              q^{12} Q^4)}{q^5 Q (q + q^3 Q) (q^5 - 
              q^6 Q^2)}E  \\
& & + \frac{(-q + q^3 Q) (q^4 + 
              q^5 Q - 2 q^6 Q - q^7 Q^2 + q^8 Q^2 - q^9 Q^2 - 
              2 q^{10} Q^3 + q^{11} Q^3 + 
              q^{12} Q^4)}{q^4 Q (q^2 + 
              q^3 Q) (-q + 
              q^6 Q^2)}E^2 \\ & &  + \frac{q^4 Q (-1 + q^3 Q)}{
(q + q^3 Q) (q - q^6 Q^2)}E^3
\end{eqnarray*}

If we wish, we may clear denominators in $P_{3_1}$ and $P_{4_1}$.
It follows that the characteristic polynomials are given by:

\begin{eqnarray*}
\ch P_{3_1}(L,M) &=&
-\frac{(L-1)(L+M^3)}{M (1 + M)} \\
\ch P_{4_1}(L,M) &=&
\frac{(L-1)(L - L M - M^2 - 2 L M^2 - L^2 M^2 - L 
M^3 + L M^4)}{M (1 + M)^2}
\end{eqnarray*}

Inspection shows that $P_{3_1}$ and $P_{4_1}$ are not regular. Nevertheless,
let us try to compute the $S$-entropy.

For the case of $3_1$, we have $|L_j(t)|=1$ for $j=1,2,3$ and in this
case
$$
\s^A_{S,3_1}(1)=\text{vol}_{3_1}=0
$$ 
for all $S$.

For $4_1$ knot, we have 3 eigenvalues $L_1(t)$, $L_2(t)=1/L_1(t)$
and $L_3(t)=1$. Assuming appropriate choices for the branches of the 
eigenvalues, the plot of $\log|L_1(t)|$ and $\log|L_2(t)|=-\log|L_1(t)|$
for $t \in [0, 2 \pi]$ is given by:  

$$
\psdraw{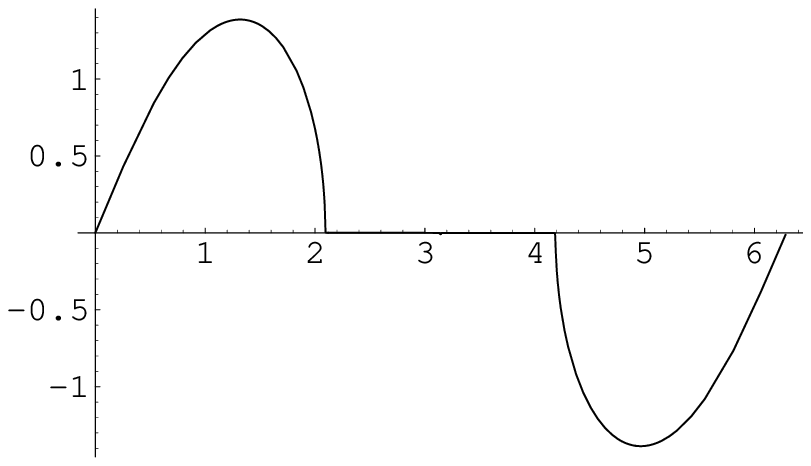}{2.7in} \qquad \psdraw{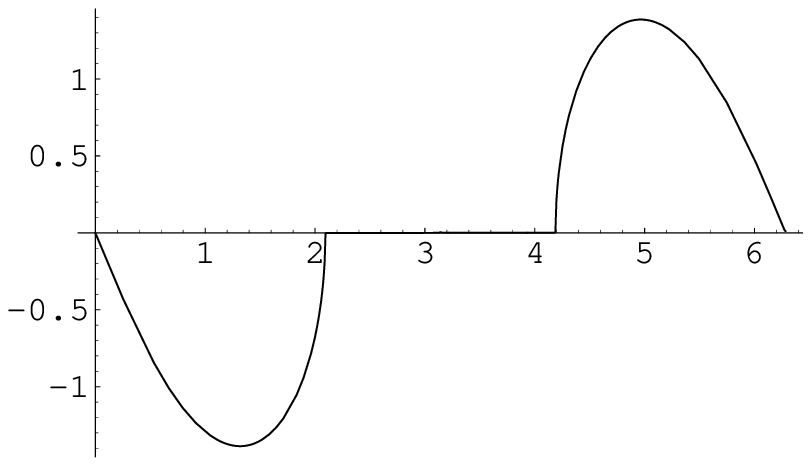}{2.7in}
$$

It follows that 
$$
\s^A_{\{1,3\},4_1}(1)=\text{vol}(4_1)=2.029883 \qquad
\s^A_{\{1,2,3\},4_1}(1) =2 \text{vol}_{4_1}=4.05977.
$$

Since the HVC is true for the $4_1$ knot, it suggests that
the colored Jones function lies in a strictly smaller subspace of the
vector space of solutions to the $q$-difference equation $P_K J_K=0$.
Using work of Murakami \cite{Mu}, one can figure out exactly the
selection principle; that is which locally fundamental solutions contribute 
to the colored Jones function.

Note that the associated $q$-difference equation of the $4_1$ knot
has the following features: of {\em collision}, {\em resonance} and 
{\em vanishing}:

\begin{itemize}
\item
The eigenvalues collide at $t=0$ (since $L_1(0)=L_1(0)=-1$), at $t=\pi/3$
(since $L_1(\pi/3)=L_1(\pi/3)=-1$) and by symmetry at $t=2 \pi/3$ and 
$t=2 \pi$.
\item
There is resonance on the interval $[\pi/3,2\pi/3]$ where all three
eigenvalues have equal magnitude.
\item
There is vanishing of the coefficients at $t=\pi/2$ 
(since the denominator $M+1$ of the coefficients is singular at $t=\pi/2$).
\end{itemize}

Moreover, there is an additional difficult problem of {\em selection 
principle}.

We plan to study these problems in later publications.

\ifx\undefined\bysame
        \newcommand{\bysame}{\leavevmode\hbox
to3em{\hrulefill}\,}
\fi

\end{document}